\documentclass[a4paper,english]{article}
\usepackage{ifthen}
\usepackage[english]{babel}
\usepackage{bibgerm}
\usepackage[utf8]{inputenc}
\usepackage{graphicx} 
\usepackage[hyphens]{url}
\usepackage{hyperref}
\usepackage{caption}
\usepackage{amsmath}
\usepackage{amsthm}
\usepackage{amsfonts}
\usepackage{amssymb}
\usepackage{mathtools}
\usepackage{enumerate}
\usepackage{xcolor}
\usepackage{xspace}
\usepackage{times}
\usepackage[a4paper, total={6in, 9in}]{geometry}
\usepackage{natbib}

\usepackage{tikz}
\usetikzlibrary{positioning}

\usetikzlibrary{shapes}
\usetikzlibrary{shapes.geometric}
\usetikzlibrary{backgrounds}
\usetikzlibrary{arrows.meta}
\usepackage{tikz-cd}
\usetikzlibrary{patterns}

\newtheorem{theorem}{Theorem}[section]

\newtheorem{definition}[theorem]{Definition}
\newtheorem{corollary}[theorem]{Corollary}
\newtheorem{remark}[theorem]{Remark}

\newtheorem{example}[theorem]{Example}
\usepackage[ruled,vlined]{algorithm2e}

\allowdisplaybreaks

\usepackage{zibtitlepage}
\usepackage{authblk}

\begin{document}

\title{Exact Objective Space Contraction for the Preprocessing of Multi-objective Integer Programs}

\author[1]{Stephanie Riedmüller\footnote{corresponding author, riedmueller@zib.de, \ZTPOrcid{0009-0006-4508-4262}}
}
\author[1,2]{Thorsten Koch
\footnote{\ZTPOrcid{0000-0002-1967-0077}}
}
\affil[1]{Zuse Institute Berlin, Berlin, Germany}
\affil[2]{Technische Universtität Berlin, Berlin, Germany}

\newcommand{\titlepageauthors}{Stephanie Riedmüller$^{1}$, Thorsten Koch$^{1,2}$}
\newcommand{\titlepageaffiliations}{
$^{1}$Zuse Institute Berlin, Germany\\
$^{2}$Technische Universtität Berlin, Berlin, Germany
}

\newcommand{\acknowledgments}{
The work for this article has been conducted in the Research Campus MODAL funded by the Federal Ministry of Research, Technology and Space (BMFTR) (fund numbers 05M14ZAM, 05M20ZBM, 05M2025).
}

\maketitle

\begin{abstract}
Solving integer optimization problems with large or widely ranged objective coefficients can lead to numerical instability and increased runtimes. When the problem also involves multiple objectives, the impact of the objective coefficients on runtimes and numerical issues multiplies. We address this issue by transforming the coefficients of linear objective functions into smaller integer coefficients. To the best of our knowledge, this problem has not been defined before. Next to a straightforward scaling heuristic, we introduce a novel exact transformation approach for the preprocessing of multi-objective binary problems. In this exact approach, the large or widely ranged integer objective coefficients are transformed into the minimal integer objective coefficients that preserve the dominance relation of the points in the objective space. The transformation problem is solved with an integer programming formulation with an exponential number of constraints. We present a cutting-plane algorithm that can efficiently handle the problem size. In a first computational study, we analyze how often and in which settings the transformation actually leads to smaller coefficients. In a second study, we evaluate how the exact transformation and a typical scaling heuristic, when used as preprocessing, affect the runtime and numerical stability of the Defining Point Algorithm.\\

\end{abstract}

\section{Introduction}


In real-world infrastructure planning and operation, Integer Linear Programs (ILPs), particularly Binary Programs with multiple objectives, become increasingly important. 

Even if a single-objective problem is easily solvable, the corresponding multi-objective problem can become increasingly hard. An attempt to classify the difficulty of multi-objective integer programs (MOIP) is given by \citet{Figueira2016} and \citet{Boekler2017}. 
In general, the difficulty increases with the number of objective functions, as investigated by \citet{Allmendinger2022}.
There is a large body of literature on algorithms for bi-objective integer programs, due to the natural ordering of efficient points along the two-dimensional Pareto front, which allows straightforward enumeration.
The challenge lies in multi-objective problems with at least three objectives.


Given a set $X \subseteq \mathbb{Z}^n$ for $n\in \mathbb{N}$, $p\geq 2$, and linear objective functions $f_1, \ldots, f_p: X \rightarrow \mathbb{R}$ with $f = (f_1, \ldots, f_p)$, we consider the multi-objective integer optimization problem 
\begin{align}
    \min_{x \in X} \quad &(f_1(x), \ldots, f_p(x)). \tag{MOIP} \label{eq:moip}\ 
\end{align}
The special case for $X \subseteq \{0,1\}^n$ is the multi-objective binary optimization problem
\begin{align}
    \min_{x \in X \subseteq \{0,1\}^n} \quad &(f_1(x), \ldots, f_p(x)). \tag{MOBO} \label{eq:moco} 
\end{align}
Recall the following definitions in the context of \ref{eq:moip}:
\begin{definition}
    We call $x \in X$ a \emph{solution} to  \ref{eq:moip} in the \emph{feasible set} $X$ of the \emph{decision space} $\mathbb{Z}^p$ and $y = f(x)$ a \emph{point} in the \emph{feasible set} $f(X)$ in the \emph{objective space} (or \emph{criterion space}) $Y$. 
\end{definition} 
\begin{definition}    
    A point $\hat{y} \in f(X)$ \emph{dominates} a point $y \in f(X)$, if $\hat{y}_i \leq y_i$ for all $i \in \{1, \ldots, p\}$ with $\hat{y}_i < y_i$ for at least one $i$. In this case, we write $\hat{y} \prec y$.
\end{definition}
Comparing elements with respect to dominance yields an irreflexive, asymmetric, and transitive binary relation, i.e., a strict partial order on the feasible set in the objective space. We refer to it as the \emph{dominance relation}.
\begin{definition}    
    A solution $\hat{x}\in X$ to  \ref{eq:moip} is called \emph{efficient} or \emph{Pareto-optimal}, if the corresponding point in the objective space is not dominated by any other point, i.e., if there is no other $x\in X$ such that $f(x) \prec f(\hat{x})$. In this case, $f(\hat{x})$ is called \emph{non-dominated point}.
\end{definition} 
\begin{definition}
    The set of efficient solutions is denoted by $X_E \subseteq X$ and is called the \emph{efficient set}. The set of non-dominated points is denoted by $Y_{N} \subseteq f(X)$ and is called the \emph{non-dominated set}.
\end{definition}

When comparing the results of multi-objective optimization algorithms, a large number of possible measures are known from the literature, most commonly the hypervolume metric. A survey of performance metrics can be found in \citet{Riquelme2015}. 

\begin{definition}
    The \emph{hypervolume (indicator)} of a set $S \subset \mathbb{R}^p$ with respect to a reference point $r \in \mathbb{R}^p$ is defined by
    \begin{equation*}
        HV(S) = \Lambda \left(\bigcup_{s\in S, r \prec s} [r_0, s_0] \times \ldots \times [r_p, s_p] \right),
    \end{equation*}
    where $\Lambda$ denotes the Lebesque measure.
\end{definition}

Hence, the hypervolume measures the union of all hyperrectangles spanned by the points in $S$ with respect to the reference point. The reference point is commonly chosen to be the ideal point.

\begin{definition}
    The point $y^I$ given by
    \begin{equation*}
        y_i^I := \min_{y\in f(X)} y_i
    \end{equation*}
    is called the \emph{ideal point}.
\end{definition}

Note that we refer to elements of the feasible set in the decision space as solutions and elements of the feasible set in the objective space as points to distinguish between the spaces.
In the following, we consider $Y := f(X)$ and then also refer to the feasible set $f(X)$ as the objective space.
For more details, we refer to 
\citet{Ehrgott2005}.

However, note that elsewhere the coefficients of the linear objective functions 
    $f_i(x) = c^\top x = \sum_{i=1}^p c_i x_i$
are ambivalently defined by $c \in \mathbb{Z}^p$ \citep{Ozlen2014, Daechert2024} or $c \in \mathbb{R}^p$ \citep{Boland2017, Bauss2023}.
While some state-of-the-art algorithms for MOIP are proven to be correct only for integer coefficients, others are theoretically correct even for continuous objective coefficients. Note that these theoretical results hold only in the idealized setting of computations with infinite precision. However, in practice, available implementations for MOIP algorithms are usually applicable only for integer coefficients.
Assuming integer coefficients without loss of generality is valid from the theoretical point of view, since continuous coefficients can be transformed into integer coefficients by scaling with a large factor. In practice, this can lead to very large and/or a wide range of objective coefficients.
The number of solutions of the corresponding problem remains the same, however:
\newpage
\begin{itemize}
    \item The hypervolume of the feasible set in the objective space increases significantly. 
    \item If the scaled coefficients have a wide range regarding the order of magnitude, numerical issues can become significant.
\end{itemize}
The same disadvantages, of course, also arise in general from large or widely ranged objective coefficients that are not generated from continuous coefficients.

In the rest of this paper, we investigate how to circumvent those issues. To that end, we recall the relevant definitions and known results about transformations in the context of multi-objective optimization. 
Essential for such a transformation is the indifference of the dominance relation under the transformation, such that the original and the transformed problem yield the same efficient set.
The following insights and the corresponding proofs can be found in \citet{NeusselStein2025}. The results are presented for linear objective functions to maintain a coherent notation throughout the paper. However, the results also hold for general continuous real-valued objective functions.

Recall that a one-dimensional function $f: \mathbb{R} \rightarrow \mathbb{R}$ is monotone if for all $x_1, x_2 \in X$, $x_1 < x_2$ implicates $f(x_1) \leq f(x_2)$ (or $f(x_1) < f(x_2)$ for strict monotonicity). This concept translates quite naturally to transformations in a multi-dimensional setting.

\begin{definition}\label{def:monotonetransformation}
    For sets $Y, Z \subseteq \mathbb{R}^p$, a mapping $\varphi: Y \rightarrow Z$ is called a \emph{monotone transformation} if it is bijective and if for all $y_1, y_2 \in Y$ holds
    \begin{align*}
        y_1 \prec y_2 \Leftrightarrow \varphi (y_1) \prec \varphi (y_2).
    \end{align*}
\end{definition}

The following result shows that monotone transformations preserve efficient solutions:
A point $\hat{y}\in Y$ is non-dominated if and only if $\varphi(\hat{y})\in Z$ is non-dominated.

\begin{theorem} \label{thm:monotone}
    For sets $X, Y, Z \subseteq \mathbb{R}^p$ and for a vector-valued linear objective function $f: X \rightarrow \mathbb{R}^p$, let $\varphi: Y \rightarrow Z$ be a monotone transformation with $f(X) \subseteq Y$. Then the sets of efficient solutions of 
        $\min_{x \in X} f(x)\text{ and of } \min_{x \in X} \varphi(f(x))$
    coincide.
\end{theorem}

\begin{definition}
    For sets $Y, Z \subseteq \mathbb{R}^p$, a bijective mapping $\varphi: Y \rightarrow Z$ is called a \emph{component-wise transformation} if it is of the form
    \begin{equation*}
        \varphi(y) = P \cdot \begin{pmatrix}
            \varphi_1(y_1)\\
            \vdots \\
            \varphi_p(y_p)
        \end{pmatrix}
    \end{equation*}
    with a permutation matrix $P$ and functions $\varphi_j: Y_j \rightarrow Z_j, j \in \{1, \ldots, p\}$. A component-wise transformation $\varphi$ is called a \emph{component-wise monotone transformation} if the functions $\varphi_j, j \in \{1, \ldots, p\}$, are monotone transformations.
\end{definition}

The following Theorem allows us to restrict the search for fitting transformations to component-wise transformations.

\begin{theorem} \label{thm:componentwise}
    For sets $Y, Z \subseteq \mathbb{R}^p$, let $\varphi: Y \rightarrow Z$ be a linear transformation. Then $\varphi$ is a monotone transformation if and only if it is a component-wise monotone transformation.
\end{theorem}

The results from Theorem~\ref{thm:monotone} and Theorem~\ref{thm:componentwise} are also true for weakly and properly efficiency, as has been proven in \citet{NeusselStein2025}. 
A solution $\hat{x}\in X$ is called \emph{weakly efficient} if there is no $x \in X$ such that $f_i(x) < f_i(\hat{x})$ for all $i \in \{1, \ldots, p\}$. In this case, $f(\hat{x})$ is called \emph{weakly non-dominated}.
There are several definitions for properly efficient and properly non-dominated. One of them is the following: A solution $\hat{x}\in X$ is called \emph{properly efficient} if it is efficient and if there is a real number $M > 0 \in \mathbb{R}$ such that for all $i$ and $x \in X$ satisfying $f_i(x) < f_i(\hat{x})$ there exists an index $j$ such that $f_j(\hat{x}) < f_j(x)$ and 
\begin{align*}
    \frac{f_i(\hat{x}) - f_i(x)}{f_j(x) - f_j(\hat{x})} \leq M.
\end{align*}
 In this case, $f(\hat{x})$ is called \emph{properly non-dominated}.
 This property allows us to consider the concepts of this paper to also apply to algorithms that provide guarantees for weakly and properly efficient solutions.

Here, we focus on preprocessing large objective spaces defined by linear objective functions.
A large objective space involves objective functions with (possibly unnecessarily) large coefficients, a wide numerical range in coefficients, or fractional objective coefficients that become large when scaled into integer values for the application of integer programming methods.
Large or widely ranged objective coefficients occur in real-world applications, for example, in industrial planning problems where high-priced investment options are combined with low-priced operational costs.
The straightforward approach to dealing with large objective coefficients is to scale and round them appropriately. 
Our contribution involves the prepossessing of large objective spaces via a novel exact transformation of the objective coefficients in multi-objective binary programs.
In contrast to the transformations in the literature, we aim not for convexification but for a reduction in the objective space hypervolume.
To the best of our knowledge, the resulting problem has not been formulated explicitly before.
We solve the transformation using an integer program with a cutting-plane algorithm.
Finally, we evaluate our method computationally. 
First, we evaluate whether our transformation can be applied to a randomly generated set of objective functions and how strongly the objectives can be contracted.
Then, we compare the application of the Defining Point Algorithm \citep{Daechert2024} to the original problem with preprocessing based on exact transformation or heuristic scaling, focusing on the number of found solutions, the number of distinct solutions, and the running time. 

In Section~\ref{section:literature}, we review the literature on preprocessing multi-objective optimization problems.
In Section~\ref{section:contraction}, we present the methodology of contracting large objective spaces via transformations, which includes a novel exact approach for the special case of binary problems.
Finally, Section~\ref{section:study} shows the impact of the presented method in two computational studies.

\section{Literature review} \label{section:literature}
 
Multi-objective integer programming is commonly divided into 
\emph{decision-space} algorithms and \emph{objective-space} algorithms. 
The decision space contains the possible alternatives, and to optimize, methods such as multi-objective branch-and-bound algorithms \citep{Bauss2023} are used.
In contrast, objective-space algorithms operate on the image of possible alternatives under the objective functions. 
The most commonly used approximation algorithms in practice are general scalarization approaches such as the weighted sum method and the $\varepsilon$-constraint method \citep{Ehrgott2006}.
Prominent examples of exact methods include approaches that divide the problem into a set of smaller subproblems and solve them using scalarization methods, see, for example, \citep{LAUMANNS2006932, Kirlik2014, Boland2017, Daechert2024}, or recursive reduction \citep{TenfeldePodehl2003, OZLEN200925, Ozlen2014}.
Among the current state-of-the-art exact objective space algorithms is the Defining Point Algorithm (DPA), presented and benchmarked by \citet{Daechert2024}, which we will use for the computations in Section~\ref{section:study}.
For an in-depth survey of exact algorithms for multi-objective integer linear optimization, we refer to \citet{Halffmann2022}, while the latest advances in exact and approximation methods are summarized by \citet{Antunes2024}.

In the following, we will present a novel preprocessing method for multi-objective binary programs. 
While preprocessing techniques for single-objective optimization are well researched \citep{Savelsbergh1994, FUGENSCHUH200569, chen2011applied} and are part of every modern solver, the literature on preprocessing for multi-objective optimization is, so far, largely restricted to two topics: dimension reduction by removing redundant/non-essential objectives \citep{GAL1977176, VAZQUEZ2018382, Kof2024, BOLAND2019858} and the convexification of non-convex Pareto sets via transformations.

The aim of dimension reduction is to identify and remove redundant objectives. The concept of redundant/non-essential objectives was introduced in \citet{GAL1977176}.
Techniques to reduce the number of objectives for general multi-objective optimization problems include 
finding a conflict error minimum objective subset (so-called $\delta$-MOSS and k-EMOSS) \citep{brockhoff2009, VAZQUEZ2018382} and
formulating an essential objective as a linear combination of the original objectives based on their correlations \citep{Cheung2014}. 
Those approaches mainly come from the community of evolutionary algorithms.
There are also methods solely for multi-objective linear programming \citep{LINDROTH20101519, Malinowska2008, Thoai2012}. A recent approach is the factorization of the coefficient matrix of low-rank objective matrices of multi-objective linear programs \citep{loehne2025}. 
For multi-objective binary problems, dimension reduction has been achieved by projections onto polyhedral cones \citep{Kof2024}.
In \citet{BOLAND2019858}, the property that a subset of the objectives of a binary optimization problem has exactly one Pareto optimal solution is exploited to combine the subset of objectives into one using a weighted sum approach.

While classical scalarization methods for multi-objective optimization transform a multi-objective problem into many single-objective problems, the literature on transforming a multi-objective problem into another multi-objective problem is rare.
Dimension reduction can also be viewed as a transformation of a multi-objective problem into one with fewer objectives.
Apart from that, transformations are used to normalize objectives with different units, especially in the weighted sum method \citep{Marler2005}.
In non-convex settings, there have been attempts to convexify Pareto fronts via transformations of the objectives involving p-powers \citep{Li1996} and exponential transformations \citep{Li1998}. 
\citet{Romeijn2004} convexifies a multi-objective problem in radiation therapy treatment planning, and \citet{Zarepisheh2017} considers convexifying transformations under strictly increasing, convex functions with a focus on the invariance of the properly efficient points. The connectedness of efficient sets of convex transformable problems is investigated by \citet{Hirschberger2005}.
Theoretical results on the component-wise structure of monotone transformations of the objective space under an invariant set of efficient points have recently been presented by \citet{NeusselStein2025} and serve as the basis for our results.
In ordinal optimization, where the standard Pareto dominance cone is substituted by a categorial ordering cone, linear transformations have been applied to the objective functions to associate the ordinal cone with the Pareto cone \citep{klamroth2023}.
Furthermore, preprocessing for multi-objective MaxSAT has been investigated by \citep{Jabs2023}

\section{Contracting large objective spaces} \label{section:contraction}

The question we are now investigating is whether it is possible to find a transformation that contracts the objective space, thereby showing that it was \emph{unnecessarily} large.
From the previous section, we know that these transformations must be bijective, monotone, and component-wise. 
Since we want to maintain the problem structure described in \ref{eq:moip} with linear objective functions, we aim at transformations that affect only the objective coefficients. The transformation should further maintain the integrality of the coefficients and, hence, of the points in the objective space to allow the application of multi-objective integer programming algorithms.

Consider the sets $X, Y, Z \subseteq \mathbb{Z}^p$ and 
\begin{align*}
    &f_i: X \rightarrow Y_i, && f_i(x) \mapsto  \sum_{j=1}^p c^i_j x_j, c^i_j \in \mathbb{N},
\end{align*} 
 the linear objective functions of \eqref{eq:moip}. The goal is to find alternative linear objective functions 
 \begin{align*}
    &g_i: X \rightarrow Z_i, && g_i(x) \mapsto  \sum_{j=1}^p d^i_j x_j, d^i_j \in \mathbb{N},
\end{align*}
such that there is a monotone, component-wise transformation $\varphi: Y \rightarrow Z$ such that $g = \varphi \circ f$, i.e., the following diagram commutes.
\[
\begin{tikzcd}
X \arrow[r, "f"] \arrow[dr, "g"'] & Y \arrow[d, "\varphi"] \\
& Z
\end{tikzcd}
\]
Here, $X$ represents the decision space, $Y$ a possibly unnecessarily large objective space, and $Z$ a possibly smaller objective space. 

When comparing the size of two objective spaces, we will use the hypervolume indicator.
So we call $Z$ smaller than $Y$ if $HV(Z) < HV(Y)$.
However, we will also consider the contraction of the individual single-objective functions.
Since the hypervolume is not applicable to single objectives, we define, as a measure for the difference in coefficients between one of the original single objective functions and the corresponding contracted single objective function, the \emph{objective contraction} by the relative change in the sum of the objective coefficients.
\begin{definition}
    Let $f,g: X \rightarrow \mathbb{Z}_{\geq 0}$ two non-negative real-valued linear functions with coefficients $c, d \in \mathbb{Z}_{\geq 0}^n$, respectively. Then the \emph{contraction factor} is defined as
    \begin{align*}
    \gamma(f,g) = \frac{\sum_{j=1}^p c_j - \sum_{j=1}^p d_j}{\sum_{j=1}^p c_j}.
\end{align*}
\end{definition}
However, not all objective functions can be contracted, as the following example shows.
\begin{example}
    The existence of such alternative non-trivial objective functions $g_j$ is not certain. Consider the biobjective problem with $X = \{0,1\}^2$ objective functions $f_1 (x) = x_1 + x_2$ and $f_2 (x) = x_1 + 2 x_2$, then there are no $g_1, g_2$ that preserve monotonicity and integrality while $Z$ is smaller than $Y$.
\end{example}
This motivates the definition of the following notion:
\begin{definition}\label{def:contractable}
    A linear integer-valued function $f: X \rightarrow \mathbb{Z}_{\geq 0}, x\mapsto \sum_{i=0}^n c_ix_i$, $c_i \in \mathbb{Z}_{\geq 0}$, $n\geq 1$ is called \emph{contractable}, if there exists another linear integer valued function $g: X \rightarrow \mathbb{Z}_{\geq 0}, x\mapsto \sum_{i=0}^n d_ix_i$, $d_i \in \mathbb{Z}_{\geq 0}$ and a monotone, transformation $\varphi: Y \rightarrow Z$ such that $g =  \varphi \circ f
    $ and $\gamma(f,g)\geq 0$. Otherwise, $f$ is called \emph{non-contractible}.
\end{definition}

Note that the monotone transformation $\varphi$ in Definition~\ref{def:contractable} must not be explicitly given in a \ref{eq:moco} setting. It is actually enough to demand that the dominance relation is preserved. 

\begin{corollary} \label{cor:implicittranformation}
    Let $X \subseteq \{0,1\}^n$, $Y, Z \in \mathbb{Z}_{\geq 0}^p$ and let $f: X \rightarrow Y$ and $g: X \rightarrow Z$ be two linear functions such that 
    for all $x_1, x_2 \in X$ holds
    \begin{align}
        f(x_1) \prec f(x_2) \Leftrightarrow  g(x_1) \prec g(x_2). \label{eq:dominance}
    \end{align}
    and
    \begin{align}
        f(x_1) = f(x_2) \Leftrightarrow  g(x_1) = g(x_2). \label{eq:welldefined}
    \end{align}
    Then the sets of efficient solutions for $x\in X$ of 
        $\min f(x)$
    and 
        $\min g(x)$ 
    coincide.
    For all $i$ with $\gamma(f_i,g_i)\geq 0$, $f_i$ is contractable.
\end{corollary}

\begin{proof}
    We assume $Y = f(X)$ and $Z = g(X)$. Define a mapping $\varphi: Y \rightarrow Z$ by $f(x) \mapsto g(x)$, which is well-defined due to \eqref{eq:welldefined}. Since $X$ is finite and therefore $Y$ and $Z$ are finite, choose  $\varphi$ to be the vector of unique polynomials $\varphi_i$ of degree $|Y| - 1$ through the distinct points $(f_i(x), g_i(x))$. Define analougously $\varphi^{-1}: Z \rightarrow Y$ by $g(x) \mapsto f(x)$. Now $\varphi$ is bijective and further a monotone transformation due to \ref{eq:dominance}. According to Theorem~\ref{thm:monotone}, the efficient sets coincide. Due to the existence of the $\varphi_i$ all $f_i$ with $\gamma(f_i,g_i)\geq 0$ are contractable.
\end{proof}

Note that $\varphi$ is indeed a component-wise monotone transformation.
In the next sections, we will make use of the property that the transformation must not be given explicitly. 
We first describe common heuristic transformations for \ref{eq:moip} and then present a new exact approach for the special case of \ref{eq:moco}.

\subsection{Heuristic approach: Scaling and rounding objective coefficients} \label{sec:scaleandround}

The brute-force approach is to scale and round the coefficients to obtain new integer coefficients. Therefore, for a given objective function $f_i(x) = \sum_{j=1}^p c_j x_j$, $c_j \in \mathbb{N}$, choose a factor $\lambda \in \mathbb{R}_{\geq 0}$ such that 
\begin{equation*}
    g_i(x) = \sum_{j=1}^p (\lambda c_j) x_j
\end{equation*}
and $d_j := (\lambda c_j) \in \mathbb{N}$. In that case, the transformation is given by 
\begin{equation*}
    \varphi: Y \rightarrow Z, \qquad y \mapsto \lambda y.
\end{equation*}
The straightforward choice to gain integrality of the new coefficients $d_j$ would be 
\begin{equation}
    \lambda = \frac{1}{gcd((c_j)_j)}. \label{eq:gcd}
\end{equation} 
However, on the one hand, the greatest common divisor tends to be $1$ for a random list of integers; on the other hand, modern solvers usually already are able to detect such an obvious scaling.

\begin{figure}
    \centering
    \scalebox{0.8}{


\begin{tikzpicture}[x=1cm,y=1cm,>=Stealth]

  \draw[step=0.5, line width=0.2pt, lightgray] (-0.2,-0.2) grid (6.7,4.2);
  \draw[step=1,   line width=0.4pt, lightgray] (-0.2,-0.2) grid (6.7,4.2);

  \draw[->] (0,0) -- (6.8,0) node[below] {$f_1$};
  \draw[->] (0,0) -- (0,4.3) node[left]  {$f_2$};

   \draw[->, teal, thick] (0,0) -- (1.5,1);
  \draw[->, violet, thick] (0,0) -- (2,2);
  \draw[->, magenta, thick] (0,0) -- (3,0);

  \foreach \P in {(2,2),(5,2), (3.5,3), (4.5,1), (6.5, 3)}
    \filldraw[black] \P circle (2.2pt);
    \foreach \Q in {(1.5,1), (3,0)}
    \draw[black] \Q circle (2.2pt);

    \node at (3,-1) {$
      f(x)
      = 
      \begin{pmatrix}
      \textcolor{violet}{4} & \textcolor{teal}{3} & \textcolor{magenta}{6} \\
      \textcolor{violet}{4} & \textcolor{teal}{2} & \textcolor{magenta}{0}
      \end{pmatrix}
      \begin{pmatrix}
        x_1 \\
        x_2 \\
        x_3 
      \end{pmatrix}
    $};

\end{tikzpicture}

\begin{tikzpicture}[x=1cm,y=1cm,>=Stealth]

  \draw[step=0.5, line width=0.2pt, lightgray] (-0.2,-0.2) grid (6.7,4.2);
  \draw[step=1,   line width=0.4pt, lightgray] (-0.2,-0.2) grid (6.7,4.2);

  \draw[->] (0,0) -- (6.8,0) node[below] {$f_1$};
  \draw[->] (0,0) -- (0,4.3) node[left]  {$f_2$};

  \draw[->, teal, thick] (0,0) -- (1.5,1);
  \draw[->, violet, thick] (0,0) -- (0.5,0.5);
  \draw[->, magenta, thick] (0,0) -- (0.5, 0);
  
  \foreach \P in {(1.5,1), (0.5, 0.5), (2, 1.5), (2,1), (1, 0.5), (2.5, 1.5)}
    \filldraw[black] \P circle (2.2pt);
    \foreach \Q in {(0.5, 0)}
    \draw[black] \Q circle (2.2pt);

    \node at (3,-1) {$
      g(x)
      = 
      \begin{pmatrix}
      \textcolor{violet}{1} & \textcolor{teal}{3} & \textcolor{magenta}{1} \\
      \textcolor{violet}{1} & \textcolor{teal}{2} & \textcolor{magenta}{0}
      \end{pmatrix}
      \begin{pmatrix}
        x_1 \\
        x_2 \\
        x_3 
      \end{pmatrix}
    $};

\end{tikzpicture}

    }
    \caption{Independent scaling of vectors in the coefficients matrix can lead to changes in the Pareto front. The feasible point set in the objective space is depicted. The non-dominated points are not filled.}
    \label{fig:scaling}
\end{figure}
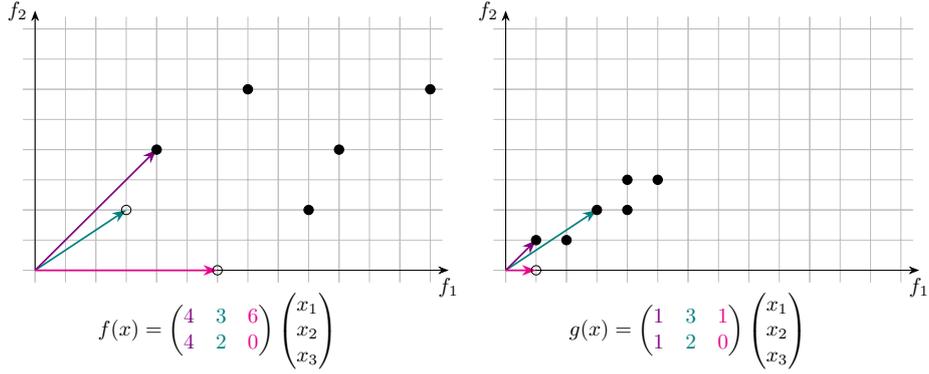

The question arises whether scaling only the row or only the column vectors is feasible.
Consider a vector-valued objective function 
\begin{align*}
    f(x) = (f_i(x))_{i = 1, \ldots, p} = \left(\sum_{j=1}^n c_{ij} x_j \right)_{i = 1, \ldots, p} = \begin{pmatrix}
        c_{11} & c_{12} & \ldots & c_{1n} \\
        c_{21} & c_{22}& \ldots & c_{2n}\\
        \vdots & \vdots & \ddots & \vdots\\
        c_{p1} & c_{p2} & \ldots & c_{pn}
    \end{pmatrix}
    \begin{pmatrix}
        x_1  \\ x_2 \\ \vdots \\ x_n
    \end{pmatrix}
    = Cx.
\end{align*}
Then the row vectors of the coefficient matrix $C$ define how the variables in the decision space are weighted in the resulting costs, i.e., they represent the optimization direction in the decision space. The column vectors represent how the corresponding single variable influences the direction in the objective space. Compare the arrows in Figure~\ref{fig:scaling} representing the column vectors. They function as a generating set of the feasible set in the objective space, since 
\begin{align*}
    Y := f(X) = \left\{ \sum_{j=0}^n c_j x_j \mid x_j \in \{0,1\}\right\},
\end{align*}
where $c_j$ denotes the $j$-th column vector. We therefore call them \emph{generating vectors}.
Scaling the column vectors of the coefficient matrix changes the lengths of the generating vectors independently. Scaling the row vectors changes the direction of the generating vectors all at once along one of the axes.

The following example shows that independently scaling column vectors can change the Pareto front. 
\begin{example}
    Consider the optimization problem
    \begin{align*}
        \min_{x \in \{0,1\}^3} \quad & \begin{pmatrix}
            4x_1 + 3 x_2 + 6x_3 \\ 4x_1 + 2x_3
        \end{pmatrix} \\
        \text{s.t.} \quad & x_1 + x_2 + x_3 \geq 1.
    \end{align*}
    The efficient set is $X_e = \{ (0,1,0), (0,0,1)\}$
    and corresponding non-dominated set $Y_n = \{  (3, 2), (6 ,0 )\}$.
    Scaling the column vectors of the coefficient matrix by $(\lambda_1, \lambda_2, \lambda_3) = (\frac{1}{4}, 1, \frac{1}{6})$, however, leads to the efficient set in $X_e = \{ (0,0,1)\}$
    and corresponding back transformed non-dominated set $Y_n = \{ (6 ,0 )\}$. See Figure~\ref{fig:scaling}.
\end{example}
In contrast, scaling the row vectors of the coefficient matrix independently is well researched. It is usually used for normalization in the weighted-sum method \citep{Marler2005}. Scaling the objective functions is, further, an inherent part of scalarization methods \citep{Ehrgott2005, Daechert2024, Helfrich2024}. In solvers for single-objective optimization, scaling the objective is usually automatically handled or can even be adjusted via a parameter. When scaling row vectors independently in a multi-objective setting, the directions of the generating vectors change synchronously along the corresponding axes; hence, the dominance relation is preserved. 

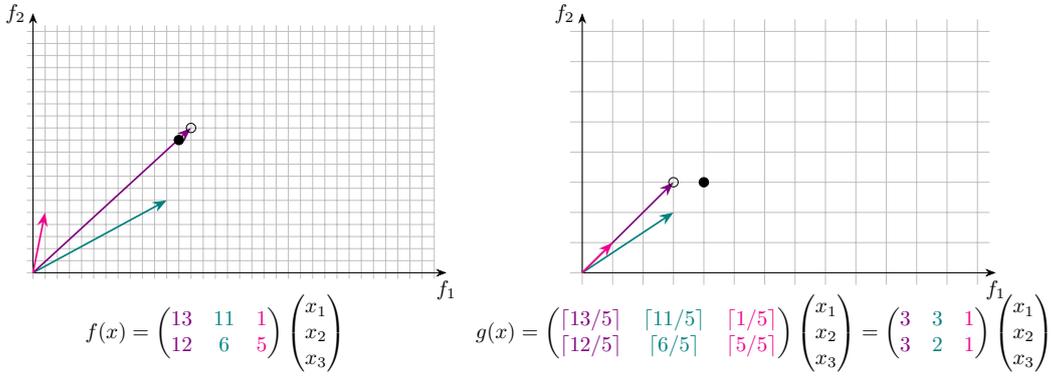
\begin{figure}
    \centering
    \scalebox{0.8}{


\begin{tikzpicture}[x=1cm,y=1cm,>=Stealth]

  \draw[step=0.2, line width=0.2pt, lightgray] (-0.1,-0.1) grid (6.6,4.1);

  \draw[->] (0,0) -- (6.8,0) node[below] {$f_1$};
  \draw[->] (0,0) -- (0,4.3) node[left]  {$f_2$};

   \draw[->, teal, thick] (0,0) -- (2.2, 1.2);
  \draw[->, violet, thick] (0,0) -- (2.6, 2.4);
  \draw[->, magenta, thick] (0,0) -- (0.2, 1);

  \foreach \P in {(2.4,2.2)}
    \filldraw[black] \P circle (2.2pt);
    \foreach \Q in {(2.6,2.4)}
    \draw[black] \Q circle (2.2pt);

  \node at (3,-1) {$
    f(x)
    =
    \begin{pmatrix}
    \textcolor{violet}{13} & \textcolor{teal}{11} & \textcolor{magenta}{1} \\
    \textcolor{violet}{12}  & \textcolor{teal}{6}  & \textcolor{magenta}{5}
    \end{pmatrix}
    \begin{pmatrix}
      x_1\\ x_2\\ x_3
    \end{pmatrix}
  $};

\end{tikzpicture}

\begin{tikzpicture}[x=1cm,y=1cm,>=Stealth]

  \draw[step=0.5, line width=0.2pt, lightgray] (-0.2,-0.2) grid (6.7,4.2);

  \draw[->] (0,0) -- (6.8,0) node[below] {$f_1$};
  \draw[->] (0,0) -- (0,4.3) node[left]  {$f_2$};

   \draw[->, teal, thick] (0,0) -- (1.5, 1);
  \draw[->, violet, thick] (0,0) -- (1.5, 1.5);
  \draw[->, magenta, thick] (0,0) -- (0.5, 0.5);

  \foreach \P in {(2,1.5)}
    \filldraw[black] \P circle (2.2pt);
    \foreach \Q in {(1.5,1.5)}
    \draw[black] \Q circle (2.2pt);

  \node at (3,-1) {$
    g(x)
    =
    \begin{pmatrix}
    \textcolor{violet}{\lceil 13/5 \rceil} & \textcolor{teal}{\lceil 11/5\rceil} & \textcolor{magenta}{\lceil 1/5\rceil} \\
    \textcolor{violet}{\lceil 12/5\rceil}  & \textcolor{teal}{\lceil 6/5\rceil}  & \textcolor{magenta}{\lceil 5/5\rceil}
    \end{pmatrix}
    \begin{pmatrix}
      x_1\\ x_2\\ x_3
    \end{pmatrix}
  =
    \begin{pmatrix}
    \textcolor{violet}{3} & \textcolor{teal}{3} & \textcolor{magenta}{1} \\
    \textcolor{violet}{3}  & \textcolor{teal}{2}  & \textcolor{magenta}{1}
    \end{pmatrix}
    \begin{pmatrix}
      x_1\\ x_2\\ x_3
    \end{pmatrix}
  $};

\end{tikzpicture}

    }
    \caption{Non-dominance is not preserved under the scale-and-round heuristic.}
    \label{fig:scaling_rounding}
\end{figure}

So let's return to a uniform scaling of the whole coefficient matrix and skip the exactness property. A straightforward heuristic is to choose $\lambda = 10^{-k}$ for some suitable $k \in \mathbb{N}$ and round the resulting $d_j$ to an integer value.
When rounding, it is recommended to round up to avoid vanishing coefficients.
This approach corresponds to scale and round fractional objective coefficients to integer values, where $d_j = \lceil (\lambda c_j) \rceil$ with $\lambda = 10^{k}$ for some suitable $k \in \mathbb{N}$.
For example, $c_j = 0.3513546$ is scaled to $3513$ using $\lambda = 10^{4}$.
The scale-and-round heuristic can diminish accuracy. 
The accuracy is diminished not only by the reduced precision of the coefficients but also by the possibility that two different coefficients map to the same value, which can affect the dominance relations of the resulting points in the objective space, as the following example shows. 
\begin{example} \label{ex:scaling_rounding}
    Consider a two-objective optimization problem with feasible set $X \subseteq \{0,1\}^3$ and objective 
    \begin{equation*}
        f(x)
    =
    \begin{pmatrix}
    13 x_1+ 11 x_2 + x_3 \\
    5 x_1  + 6 x_2  + 10 x_3
    \end{pmatrix}.
    \end{equation*}
    For solution vectors $x=(0,1,1)$ and $y=(1,0,0)$, the corresponding points yield $f(x) \prec f(y)$. However, when $f$ is transformed by the scale-and-round heuristic into $g$ by a scaling by $5$, i.e. $g(x) = \sum_{i=1}^3 \lceil \frac{c_i}{5} \rceil x_i$, the relation changes to $g(y) \prec g(x)$. Compare Figure~\ref{fig:scaling_rounding}.
\end{example}

Furthermore, it is not straightforward to determine which accuracy is required for a given case.
The loss of accuracy when transforming continuous values into integers of limited precision is well studied in numerical analysis under the term of discretization error and in signal processing under the term of quantization error \citep{Higham2002, Widrow_Kollár_2008}.

\subsection{Exact approach: Mapping objective coefficients in binary programs} \label{sec:contraction}

\begin{figure}
    \centering
    \scalebox{0.8}{


\begin{tikzpicture}[x=1cm,y=1cm,>=Stealth]

    \filldraw[fill=gray!30, fill opacity=0.5, draw=gray!30]
    (0,0) -- (0.5, 0) -- (0.5,2.5)  -- (0, 2.5) --  cycle;
    \filldraw[fill=gray!30, fill opacity=0.45, draw=gray!30]
    (0,0) -- (1, 0) -- (1, 0.5)  -- (0, 0.5) --  cycle;
    \filldraw[fill=gray!30, fill opacity=0.4, draw=gray!30]
    (0,0) -- (5.5, 0) -- (5.5, 0.5)  -- (0, 0.5) --  cycle;
    \filldraw[fill=gray!30, fill opacity=0.35, draw=gray!30]
    (0,0) -- (5, 0) -- (5, 2.5)  -- (0, 2.5) --  cycle;
    \filldraw[fill=gray!30, fill opacity=0.3, draw=gray!30]
    (0,0) -- (1.5, 0) -- (1.5, 3)  -- (0, 3) --  cycle;
    \filldraw[fill=gray!30, fill opacity=0.25, draw=gray!30]
    (0,0) -- (6, 0) -- (6, 3)  -- (0, 3) --  cycle;
    
  \draw[step=0.5, line width=0.2pt, lightgray] (-0.2,-0.2) grid (6.2,4.2);
  \draw[step=1,   line width=0.4pt, lightgray] (-0.2,-0.2) grid (6.2,4.2);

  \draw[->] (0,0) -- (6.3,0) node[below] {$f_1$};
  \draw[->] (0,0) -- (0,4.3) node[left]  {$f_2$};

   \draw[->, teal, thick] (0,0) -- (0.5,2.5);
  \draw[->, violet, thick] (0,0) -- (1, 0.5);
  \draw[->, magenta, thick] (0,0) -- (4.5, 0);

  \foreach \P in {(0,0),(0.5,2.5), (1, 0.5),(4.5, 0), (5.5, 0.5),(5, 2.5),(1.5, 3),(6, 3)}
    \filldraw[black] \P circle (2.2pt);

    \node at (3,-1) {$
      f(x)
      = 
      \begin{pmatrix}
      \textcolor{violet}{2} & \textcolor{teal}{1} & \textcolor{magenta}{9} \\
      \textcolor{violet}{1} & \textcolor{teal}{5} & \textcolor{magenta}{0}
      \end{pmatrix}
      \begin{pmatrix}
        x_1 \\
        x_2 \\
        x_3 
      \end{pmatrix}
    $};

    \fill[pattern=north east lines, pattern color=gray]
    (1.8,-0.2) rectangle (4.2,4.2);
    \fill[pattern=north east lines, pattern color=gray]
    (-0.2,0.8) rectangle (6.2,2.2);

\end{tikzpicture}

\begin{tikzpicture}[x=1cm,y=1cm,>=Stealth]

    \filldraw[fill=gray!30, fill opacity=0.5, draw=gray!30]
    (0,0) -- (0.5, 0) -- (0.5,1.5)  -- (0, 1.5) --  cycle;
    \filldraw[fill=gray!30, fill opacity=0.45, draw=gray!30]
    (0,0) -- (1, 0) -- (1, 0.5)  -- (0, 0.5) --  cycle;
    \filldraw[fill=gray!30, fill opacity=0.4, draw=gray!30]
    (0,0) -- (2.5, 0) -- (2.5, 0.5)  -- (0, 0.5) --  cycle;
    \filldraw[fill=gray!30, fill opacity=0.35, draw=gray!30]
    (0,0) -- (2, 0) -- (2,1.5)  -- (0, 1.5) --  cycle;
    \filldraw[fill=gray!30, fill opacity=0.3, draw=gray!30]
    (0,0) -- (1.5, 0) -- (1.5,2)  -- (0, 2) --  cycle;
    \filldraw[fill=gray!30, fill opacity=0.25, draw=gray!30]
    (0,0) -- (3, 0) -- (3,2)  -- (0, 2) --  cycle;

  \draw[step=0.5, line width=0.2pt, lightgray] (-0.2,-0.2) grid (6.2,4.2);
  \draw[step=1,   line width=0.4pt, lightgray] (-0.2,-0.2) grid (6.2,4.2);
    
  \draw[->] (0,0) -- (6.3,0) node[below] {$f_1$};
  \draw[->] (0,0) -- (0,4.3) node[left]  {$f_2$};

  \draw[->, teal, thick] (0,0) -- (0.5,1.5);
  \draw[->, violet, thick] (0,0) -- (1, 0.5);
  \draw[->, magenta, thick] (0,0) -- (1.5, 0);
  
  \foreach \P in {(0,0),(0.5,1.5), (1, 0.5),(1.5, 0), (2.5, 0.5),(2,1.5),(1.5,2),(3, 2)}
    \filldraw[black] \P circle (2.2pt);

    \node at (3,-1) {$
      g(x)
      = 
      \begin{pmatrix}
      \textcolor{violet}{2} & \textcolor{teal}{1} & \textcolor{magenta}{3} \\
      \textcolor{violet}{1} & \textcolor{teal}{3} & \textcolor{magenta}{0}
      \end{pmatrix}
      \begin{pmatrix}
        x_1 \\
        x_2 \\
        x_3 
      \end{pmatrix}
    $};

\end{tikzpicture}

    }
    \caption{Contraction of an unnecessarily large objective space for a binary linear problem with two objectives. The filled area represents the hypervolume of the feasible set in the objective space, and the hatched area represents the contracted space.}
    \label{fig:contraction}
\end{figure}
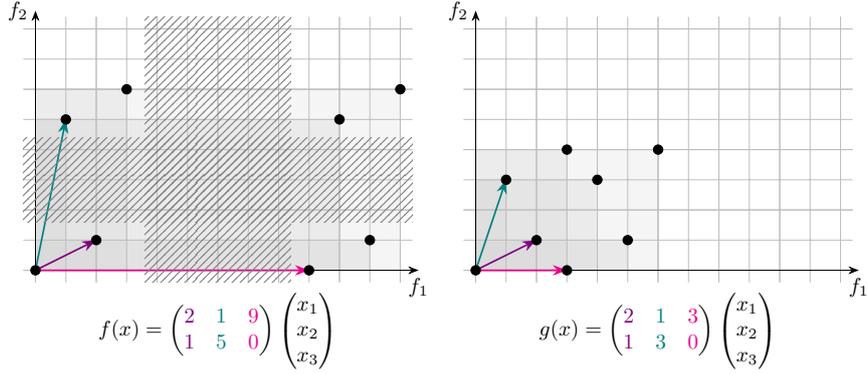

When dealing with an integer program, all integer points in the feasible space of the linear programming relaxation are actually feasible. However, when dealing with a binary program, the feasible space becomes necessarily finite and more sparse in the sense that integer points in the feasible space of the linear programming relaxation are no longer necessarily feasible. 
We suggest that those infeasible regions can be determined by the objective coefficients and, moreover, can be cut out.
Figure~\ref{fig:contraction} illustrates the problem for a binary linear program with two objectives in three variables. The arrows indicate the generating vectors that span the feasible objective space given by the components of the two objective functions $f = (f_1, f_2)^\top$. 
However, the objective space spanned by $g = (g_1, g_2)^\top$ has an equivalent structure, in the sense that the dominance relation between the feasible points does not change. When using $g$ instead of $f$, the gray hatched area, which includes no feasible points, is contracted.

In the following, without loss of generality, we consider all coefficients to be positive integers sorted in increasing order.

\subsubsection{Full model}

We need a transformation that maps the coefficients $c \in \mathbb{Z}_{\geq 0}^n$ of the objective function to other coefficients $d \in \mathbb{Z}_{\geq 0}^n$ while preserving non-dominance. Since the transformation must be component-wise according to Theorem~\ref{thm:componentwise}, we can consider each objective function separately. Note that this contraction is, hence, also valid for single-objective problems. According to Corollary~\ref{cor:implicittranformation}, the required transformation must not be explicitly given. Hence, it is enough if we can determine the transformed objective function by defining its coefficients $d \in \mathbb{Z}_{\geq 0}^n$. To ensure non-dominance preservation for the transformation of the multi-objective function, the transformed individual objectives must preserve monotonicity in each component. Monotonicity is preserved if the following implications hold for binary vectors $x,y$: if $f(x) < f(y)$ then $g(x) < g(y)$ and if $f(x) = f(y)$ then $g(x) = g(y)$. Note that the equality implication is only necessary if $f$ is not injective. 

Given $X \subseteq \{0,1\}^n$ and $c \in \mathbb{Z}_{\geq 0}^n$, we define variables $d_i \in \mathbb{Z}_{\geq 0}, i = 1, \ldots , n$. The \emph{objective contraction problem (OCP)} of finding a minimal $d \in \mathbb{Z}_{\geq 0}^n$ can be formulated by the following integer program:
\begin{align*}
    \min \quad & \sum_{i=0}^n d_i \tag{OCP} \label{eq:ocp}\\
    &d^{\top}x - d^{\top}y \leq - 1&&\forall (x,y) \in X^2: c^{\top}x < c^{\top}y \\
    &d^{\top}x - d^{\top}y = 0&&\forall (x,y) \in X^2: c^{\top}x = c^{\top}y \\
    &d_i \in [1, \max_i c_i] \cap \mathbb{Z} &&\forall i = 1, \ldots , n
\end{align*}
This formulation includes only $n$ variables. The number of constraints is exponential, since all possible combinations of solution vectors $(x,y) \in X^2$ need to be enumerated.
\ref{eq:ocp} can be reformulated based on sets:
\newpage
\begin{align}
    \min \quad & \sum_{i=0}^n d_i  \tag{OCPset}\label{eq:setocp}\\
    &\sum_{i \in S} d_i - \sum_{i \in T} d_i  \leq 1&&\forall S, T \subseteq [n], S\cap T = \emptyset, \text{ if } \sum_{i \in S} c_i < \sum_{i \in T} c_i \\
    &\sum_{i \in S} d_i - \sum_{i \in T} d_i = 0 &&\forall S, T \subseteq [n], S\cap T = \emptyset, \text{ if } \sum_{i \in S} c_i = \sum_{i \in T} c_i \label{eq:setequality}\\
    &d_i \in [1, \max_i c_i] \cap \mathbb{Z} &&\forall i = 1, \ldots , n
\end{align}

\begin{theorem}
    For given $c \in \mathbb{Z}_{\geq 0}$, let $P_{\ref{eq:ocp}}$ be the set of feasible solutions to \eqref{eq:ocp} and $P_{\ref{eq:setocp}}$ be the set of feasible solutions to \eqref{eq:setocp}. Then $P_{\ref{eq:ocp}} = P_{\ref{eq:setocp}}$.
\end{theorem}

\begin{proof}
    $d \in P_{\ref{eq:ocp}}$ if and only if $d \in P_{\ref{eq:setocp}}$ follows by case distinction and choosing $S = \{i \mid x_i = 0, y_i = 1 \}$ and $T = \{i \mid x_i = 1, y_i = 0 \}$ appropriately.    
\end{proof}

\subsubsection{Cutting-plane algorithm}

We now present a cutting-plane algorithm to solve \ref{eq:setocp}.
The algorithm is described in pseudocode in Algorithm~\ref{alg:cuttingplain}. 
Due to the separation procedure, only those of the exponential many constraints that are violated are added to the linear programming relaxation. This leads to a substantial speed-up in practice. 

For the relaxed problem, include only the constraints for the relation of subsequent coefficients, since those are directly accessible from the input coefficients. Recall that the $c_i$ are sorted in ascending order.
Depending on the relation of $c_i$ and $c_{i+1}$, add 
\begin{equation*}
    d_i + 1 \leq d_{i+1} \qquad \text{ or } \qquad d_i = d_{i+1},
\end{equation*} 
for all $i=0, \ldots, n$.
Now we iteratively add those constraints that are (maximally) violated.
Violated constraints can be found by the following separation oracle (see Algorithm~\ref{alg:oracleA}):
The task is to find sets of indices $S, T, S\cap T = \emptyset$ such that $\sum_{i \in S} d_i \geq \sum_{i \in T} d_i $ and $\sum_{i \in S} c_i \leq \sum_{i \in T} c_i $.
To that end, define variables
\begin{align*}
    z_i = 
    \begin{cases}
    1 & \text{if } i \in S \\
    -1 & \text{if } i \in T \\
    0 & \text{otherwise }
\end{cases}
\end{align*}
and solve the following IP:
\begin{align*}
    \max \quad & \sum_{i=0}^n d_i z_i \tag{OrA} \label{eq:orA}\\
    s.t. \quad&\sum_{i=0}^n c_i z_i  \leq 0\\
    &z_i \in \{-1, 0 ,1\}
    \end{align*}
If the objective value is positive, the index sets $S$ and $T$ can be reconstructed such that we generate a violated constraint 
\begin{align*}
    \sum_{i \in S} d_i +1 \leq \sum_{i \in T} d_i &&\text{or } &&&\sum_{i \in S} d_i = \sum_{i \in T} d_i
\end{align*}
depending on the value of $\sum_{i=0}^n c_i z_i$.
For the case of a vanishing objective value, there is no constraint of the previous form violated. However, \ref{eq:orA} misses a critical case, when the objective vanishes; hence, we need to consider a second stage separation oracle (see Algorithm~\ref{alg:oracleB}):
\newpage
\begin{align*}
    \min \quad & \sum_{i=0}^n c_i z_i \tag{OrB} \label{eq:orB}\\
    s.t. \quad&\sum_{i=0}^n d_i z_i  = 0\\
    &z_i \in \{-1, 0 ,1\}
    \end{align*}
For a negative objective, we generate the cut
\begin{align*}
    \sum_{i \in S} d_i +1 \leq \sum_{i \in T} d_i.
\end{align*}
The stopping criteria are vanishing objective values for both oracles \ref{eq:orA} and \ref{eq:orB} or infeasibility.

Note that the first stage separation oracle \eqref{eq:orA} is essentially a Knapsack problem, and the second stage separation oracle \eqref{eq:orB} is essentially a Subsetsum problem. While both problem classes are known to be NP-hard, the size of the oracle problems is manageable.

    The described separation oracles involve integer variables. A corresponding binary formulation can be constructed by introducing binary variables $x_i, y_i \in \{0,1\}$ with $x_i + y_i \le 1$ and substituting $z_i = x_i - y_i$ for all $i=0,\dots,n$, such that 
\begin{align*}
\begin{aligned}
y_i &= 
\begin{cases}
  1 & \text{if } i \in S \\
  0 & \text{otherwise }
\end{cases}
&
\quad \text{ and } \quad
z_i &=
\begin{cases}
  1  & \text{if } i \in T \\
  0 & \text{otherwise. }
\end{cases}
\end{aligned}
\end{align*}
Note that the variables $x,y$ directly represent the solution vectors with the same name in \ref{eq:ocp}. In the following, we use the integer version.

\begin{algorithm}[ht] 
\caption{Objective function contraction}
\label{alg:cuttingplain}
\KwIn{positive integer vector $c$ with increasing entries, relaxed problem $R$}
\KwOut{minimal positive integer vector $d$}
\ForEach{$i = 0, \ldots n-1$}{
    \eIf{$c_i = c_{i+1}$}{
        add constraint $d_i = d_{i+1}$ to $R$;
    }{
        add constraint $d_i + 1 \leq  d_{i+1}$ to $R$;
    }
}
\Repeat {$\mu =$ None an $\nu =$ None} {
$d = solve (R)$\;
$\mu =$ SeparationOracleA($d$)\;
    \eIf{$\mu$ not None}{
        Add cut $\mu$ to $R$\;
    }{
        $\nu =$ SeparationOracleB($d$)\;
        \If{$\nu$ not None}{
        Add cut $\nu$ to $R$\;
    }
    }
    }
\Return None\;
\end{algorithm}

\begin{algorithm}[ht] 
\caption{Separation Oracle A}
\label{alg:oracleA}
\KwIn{temporal solution $d$}
\KwOut{violated inequality or None}
$\delta = \max \sum_{i=0}^n d_i x_i$
    s.t. $\sum_{i=0}^n c_i x_i  \leq 0, x_i \in \{-1, 0 ,1\}$\;
\If{$\delta > 0$}{
$\kappa = \sum_{i=0}^n c_i x_i $\;
$S = \{i \mid x_i = 1 \}$, $T= \{i \mid x_i = -1 \}$\;
\If{$S, T \ne \emptyset$}{
        \eIf{$\kappa = 0$}{
        \Return $\sum_{i \in S} d_i = \sum_{i \in T} d_i$\;
    }{
        \Return $\sum_{i \in S} d_i + 1\leq \sum_{i \in T} d_i$\;
    }
    }
    }
\Return None\;
\end{algorithm}

\begin{algorithm}[ht] 
\caption{Separation Oracle B}
\label{alg:oracleB}
\KwIn{temporal solution $d$}
\KwOut{violated inequality or None}
$\kappa = \min \sum_{i=0}^n c_i x_i$
    s.t. $\sum_{i=0}^n d_i x_i  = 0, x_i \in \{-1, 0 ,1\}$\;
\If{$\kappa < 0$}{
    $S = \{i \mid x_i = 1 \}$, $T= \{i \mid x_i = -1 \}$\;
        \Return $\sum_{i \in S} d_i + 1\leq \sum_{i \in T} d_i$\;
    }
\Return None\;
\end{algorithm}

\subsubsection{Existence of contractable objective functions }

The following Theorem shows that, independent of the number of variables, there always exist contractable and non-contractable objective functions. 

\begin{theorem} \label{thm:existence}
    For all $n \in \mathbb{Z}_{\geq 0}$,  
    \begin{enumerate}
        \item there exists a vector $c \in \mathbb{Z}_{\geq 0}^n$ such that $d=c$ is not optimal for \ref{eq:setocp}, i.e., a linear function $f$ with coefficients $c$ is contractable.
        \item there exists a vector $c \in \mathbb{Z}_{\geq 0}^n$ such that $d=c$ is optimal for \ref{eq:setocp}, i.e., a linear function $f$ with coefficients $c$ is non-contractable.
    \end{enumerate}
\end{theorem}

\begin{proof}
    Let $n \in \mathbb{Z}_{\geq 0}$ be arbitrary. For both parts of the theorem, we construct a sequence $(c_i)_{i=1,\ldots,n}$ that fulfills the statement. Since $n \leq 2 $ is trivial, consider $n\geq 3$:
    \begin{enumerate}
        \item Consider $c_i = t \cdot 2^i$ for some $t \in \mathbb{Z}, t \geq 2$. Let $S, T \subseteq [n], S\cap T = \emptyset$. If $$\sum_{i \in S} c_i  = \sum_{i \in S} t \cdot 2^i < \sum_{i \in T} t \cdot 2^i  = \sum_{i \in T} c_i$$ holds, then also $$\sum_{i \in S} d_i  = \sum_{i \in S} 2^i < \sum_{i \in T} 2^i  = \sum_{i \in T} d_i.$$ The same argument holds for equality.
        \item Consider $c_i = 2i + 1$. Let $S = \{k, k+ t+ 3\}$ and $T = \{k+1, k + t +2\}$ for $k,t \in \mathbb{Z}_{\geq 0}$, then $$\sum_{i \in S} c_i  = 4k +2t + 8 =  \sum_{i \in T} c_i.$$ Using Constraint~\ref{eq:setequality} leads to 
        $$d_{k+1} - d_k  = d_{k+t+3} - d_{k+t+2}$$
        and for $k=0$ to
        $$d_{1} - d_0  = d_{t+3} - d_{t+2}.$$
        Hence, all succeeding coefficients $d_i$ must have the same difference $\delta_d = d_1 - d_0$. Since $d_0 < d_1 \in \mathbb{Z}_{\geq 0}$, $\delta_c = c_1 - c_0 = 2$ and since the sum of $d_i$ is minimized, this yields $\delta_d \in \{1,2\}$. For $S= \{0,2 \}, T = \{ 3\}$, the constraint $1 \leq d_0 < d_3 - d_2 = \delta_d$, the difference must be $\delta_d = 2$. Hence, $d=c$ is optimal.
    \end{enumerate}
\end{proof}

Even though for each $n$ a linear function with coefficient vector $c$ exists that is contractable, the probability of an objective function being contractable decreases with increasing $n$. Given a problem with $n$ variables, an additional variable can lead to exponentially many new constraints in problem \ref{eq:setocp}, since the number of possible subsets $S \in [n]$ increases by $2^n$.

\begin{remark}
    A short excursus: In the theory of positional numeral systems, as the binary or decimal system, the radix or basis $b$ is the number of unique digits per position. Typically, in such a system the string of digits $d_1 ... d_n$ denotes the number $d_1b^{n-1} + d_2b^{n-2} + ... + d_nb^0$, where $0 \leq d_i < b$. If the basis is fixed, the numeral system is called \emph{standard} and for a mixed basis \emph{non-standard}. A standard numeral system fulfills the property $\sum_{i< k}d_i b^{i} < b^k$ for all $k$ and therefore allows for a unique representation of each number.    
    The notion of (non-)standard numeral systems can be transferred to our objective contraction problem in the following way: The objective coefficients $c_i$ represent the values of the positions $b^i$, and the solution vector components $x_i$ represent the digits $d_i$. A binary vector $x$ can then be interpreted as the string of digits denoting the objective value $\sum c_i x_i$. Since the coefficients are usually not identical, this would result in a non-standard numeral system. 
    However, if we consider a coefficient vector with $\sum_{i < k} c_i < c_k$ for all $k$, then this vector also has the property of unique representation as a standard numeral system. We, therefore, call it a \emph{standard numeral} coefficient vector with standard numeral coefficients.
    A vector of coefficients induces an order on the objective values. If the coefficient vector is a standard numeral, then the induced order is the inverse lexicographic order, i.e., standard coefficients and strict lexicographic dominance are equivalent.
    To see this, let $x <_{lex} y\in \{0,1\}^n$ with $x_i = y_i = 0$ for $i>k$ and let $c_k$ be standard numeral. Then also $\sum c_i x_i < \sum c_i y_i$, since $c_k$ dominates the sum by $\sum_{i < k} c_i < c_k$. On the other hand, if $c_k$ is not standard numeral, then $\sum_{i < k} c_i \geq c_k$ and therefore strict lexicographic dominance does not hold for example for $x$ with $x_{i < k} = 1, x_{i \geq k} = 0$ $y$ with $y_k = 1, y_{i \ne k} = 0$.
    The coefficient sequence $d_i=2^i$ from Theorem~\ref{thm:existence} is, therefore, optimal and also the optimal solution for all standard numeral coefficient vectors.
\end{remark}

\subsubsection{Negative coefficients}

In the previous sections, we assumed the coefficients to be positive integers. For negative coefficients, the approach is still valid but requires some minor adjustments to problem \ref{eq:setocp}. Minimizing the sum over all new objective coefficients $d_i, i=0, \ldots, n$ does not minimize the objective space, since large negative coefficients could compensate large positive coefficients. This leads to two possible alternatives:
\begin{enumerate}
    \item Square the individual coefficients, i.e., minimize $\sum_{i=0}^n d_i^2$. This, however, results in a quadratic problem that must either be solved with a nonlinear solver or linearized, which introduces more variables.
    \item Divide coefficients in positive and negative part $d_i = p_i - n_i$, $p_i, n_i \in \mathbb{Z}_{\geq 0}$. Then minimize $\sum_{i=0}^n p_i + n_i$. This approach results in twice as many variables but avoids destroying linearity.
\end{enumerate}
In both cases, the separation oracles remain unaffected.

\section{Computational study} \label{section:study}

For the following studies, all computations have been carried out on an Intel Xeon CPU E3-1245 v6 @ 3.70GHz using C++ and the CPLEX solver version 12.10.0 \citep{ibm_cplex}.
The cutting-plane algorithm has been implemented with lazy callbacks and executed on a single CPU to ensure callback safety. The trivial solution $d = c$ is given to the solver as a warm-start.

\subsection{Contraction of individual objectives}

For the investigation of the impact of contractions on objective coefficients, we first restrict to single objectives. The goal is to determine how often objectives are indeed contractable. We generate positive integer objective coefficients randomly with three different variants:
\begin{enumerate}
    \item Uniform: The coefficients are integers drawn uniformly at random in a given range. A uniform sampling from a large range, however, has a bias towards large numbers.
    \item OrderOfMagnitude: Precompute all orders of magnitude that intersect the target range. Then, for each coefficient, randomly select one of these orders of magnitude and sample a value uniformly from within that order. This approach produces coefficients whose magnitudes are evenly distributed across the orders of magnitude.
    \item LogarithmicSpread: Sample coefficients logarithmically over the range such that there is a natural bias towards smaller numbers within a large range.
\end{enumerate}

For a number of coefficients between 5 and 30, and for a range $(1, 10^k) $ for $k = 3, \ldots 7$, we generate 5 samples from each of the described variants of sampling. This results in 75 samples per coefficient. A time limit of 600s is set.

\begin{figure}
    \centering
    \includegraphics[width=\linewidth]{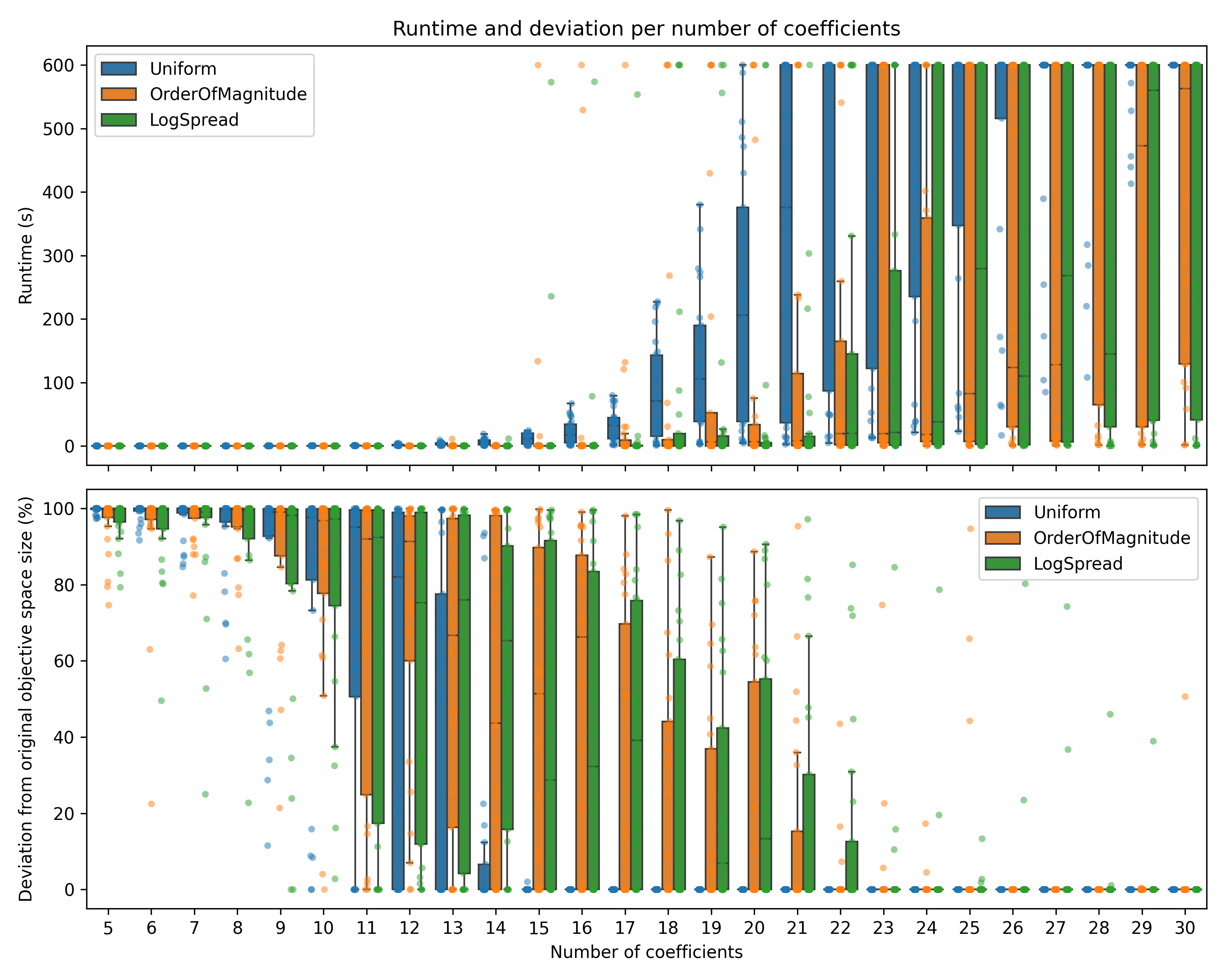}
    \caption{Runtime and size deviation (contraction factor multiplied by 100) of the contraction of objective coefficients depicted as a box plot per number of coefficients.}
    \label{fig:so_boxplot}
\end{figure}

Figure~\ref{fig:so_boxplot} depicts the runtime of the contraction as well as the deviation in the sum of the objective coefficients (contraction factor $\gamma(f,g)$ given in percentage). The color indicates the sampling strategy. 
The impact of the contraction on the size of coefficients is especially high for up to $10$ coefficients (median lies above 95\% for all $n\leq 10$) but decreases with increasing $n$. 
From $23$ coefficients, the impact of the contraction is only visible in outliers. The uniform distribution of coefficients shows a consistent impact up to $14$ coefficients.
Hence, the chance that the contraction can be applied successfully increases if the coefficients are not only large but also wide in range.
The runtime shows an opposite development as expected due to the increase in problem size. 

\begin{figure}
    \centering
    \includegraphics[width=0.8\linewidth]{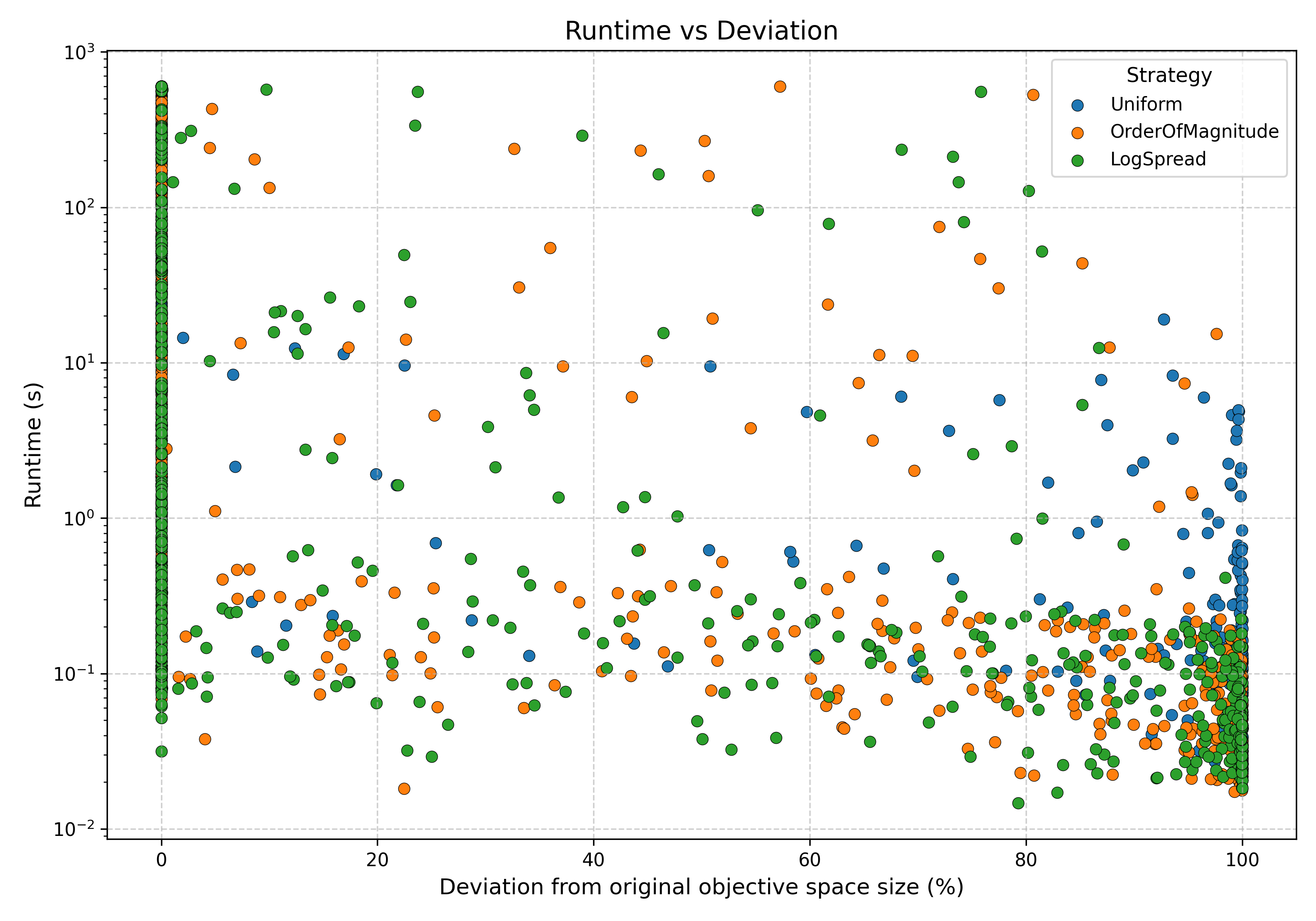}
    \caption{Runtime vs deviation in coefficient size (contraction factor multiplied with 100) for a single objective contraction.}
    \label{fig:so_runtime_deviation}
\end{figure}
\begin{figure}
    \centering
    \includegraphics[width=0.8\linewidth]{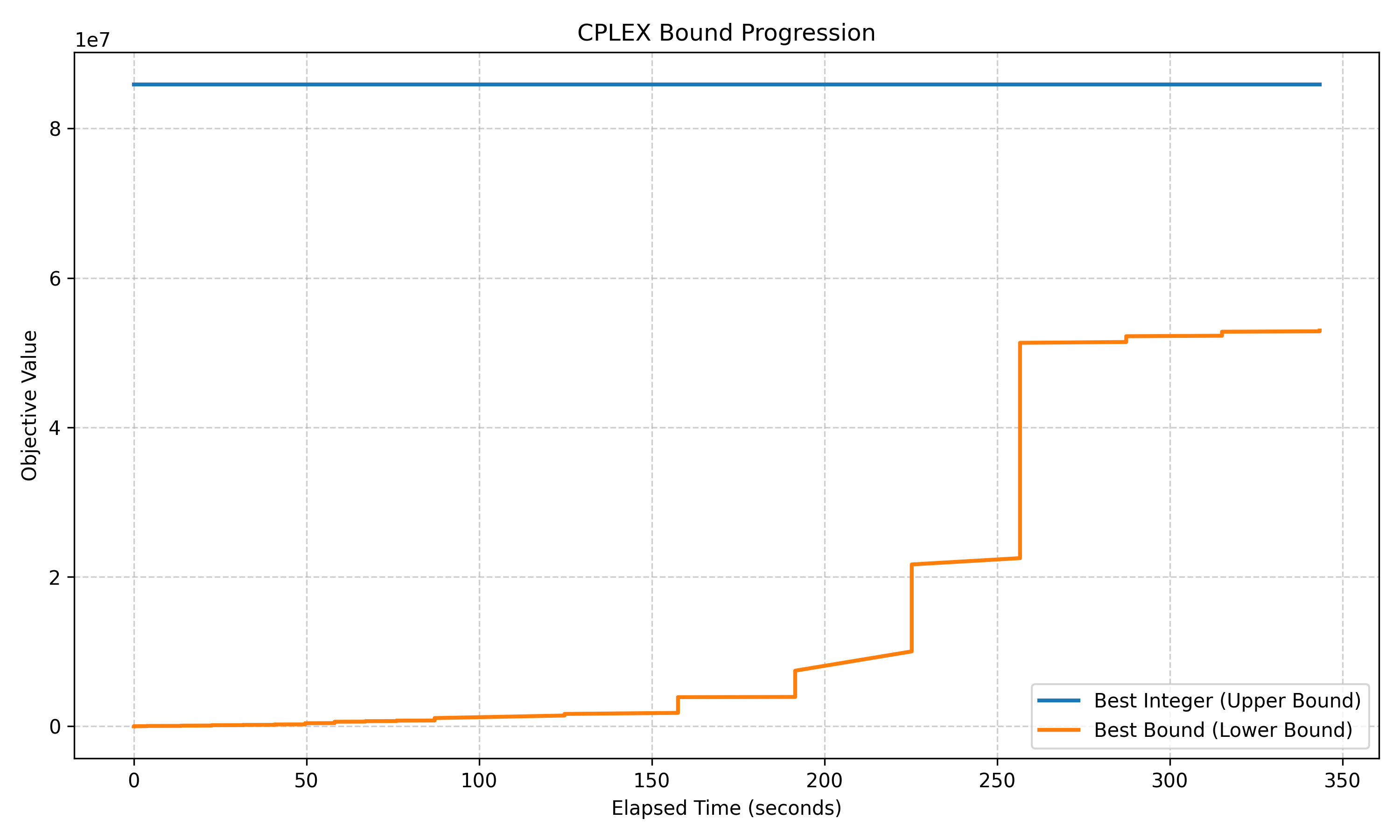}
    \caption{Development of the values for the lower and upper bound for a non-contractable objective function.}
    \label{fig:bounds}
\end{figure}

Figure~\ref{fig:so_runtime_deviation} shows the runtime dependent on the deviation (contraction factor given in percentage) and underlines the observation that high runtimes appear mainly for non-contractible objective functions. If problem \ref{eq:setocp} actually has an optimal solution $d* \ne c$, the solution process using the cutting-plane algorithm solves efficiently. However, if the optimal solution is $d^* = c$, i.e., is identical to the warm-start solution, the solving time is strongly increased. 
Hence, finding the optimal objective coefficients is easy, while proving optimality is the bottleneck.
Figure~\ref{fig:bounds} shows an exemplary plot of the upper and lower bounds for a non-contractible objective function, i.e., from a case where the previous objective already was optimal. The solver then stops because all nodes of the branch-and-bound tree have been investigated.

\subsection{Contraction of the objective space}

For the following computational study, we used the openly available implementation of DPA \citep{Daechert2024}. We consider DPA as a representative of state-of-the-art objective space integer optimization algorithms. DPA is an exact objective space algorithm that iteratively decomposes the objective space into non-disjoint rectangular boxes. The boxes are described by upper-bound vectors. The decomposition is redefined in each iteration depending on the previously found non-dominated point.
Each subproblem, i.e., each box, is solved with a scalarization method by calling the commercial MIP solver CPLEX version 12.10.0 \citep{ibm_cplex}. We choose as scalarization the augmented $\varepsilon$-constraint method (DPA-a), which combines the multi-objective subproblem into a single objective problem, by restricting the objective functions through the box constraints and optimizing for
\begin{equation*}
    f_0(x) + \rho \sum_{j=1}^m f_j(x)
\end{equation*}
with a small augmentation parameter $\rho > 0$.
The algorithm guarantees a minimal number of local upper bound vectors and no redundancy in the generation of non-dominated points.

\subsubsection{Instability in DPA}

For large or widely ranged coefficients, DPA tends to regenerate duplicate solutions unnecessarily. As a result, the solution process can even become stuck and fail to terminate within a reasonable time. We briefly illustrate how this occurs with the following example.

\begin{figure}
    \centering


\begin{tikzpicture}[x=1cm,y=1cm,>=Stealth]

  \draw[step=0.5, line width=0.2pt, lightgray] (-0.2,-0.2) grid (6.7,4.2);

  \draw[->] (0,0) -- (6.8,0) node[below] {$f_1$};
  \draw[->] (0,0) -- (0,4.3) node[left]  {$f_2$};

  \fill[pattern=north east lines, pattern color=gray]
    (0.5, 0.5) rectangle (6.5, 2.75);

   \draw[-, teal, thick] (-0.2,0.5) -- (6.7, 0.5);
   \draw[-, teal, thick] (-0.2,4) -- (6.7, 4);
   \draw[-, teal, thick] (6.5,4.2) -- (6.5, -0.2);
   \draw[-, teal, thick] (0.5,-0.2) -- (0.5, 4.2);

   \draw[dotted, magenta, thick] (3.75,-0.2) -- (3.75,4.2);
   \draw[dotted, magenta, thick] (-0.2,2.75) -- (6.7,2.75);

    \filldraw[black] (4,3) circle (2.2pt) node[above right] {$f(x)$};
    \draw[black] (5,2.5) circle (2.2pt) node[below right] {$f(y)$};

\end{tikzpicture}

    \caption{Problematic regeneration of redundant solutions in DPA. The solid lines define the search area, the dotted lines define the boxes, and the hatched area represents the current box.}
    \label{fig:dpa}
\end{figure}
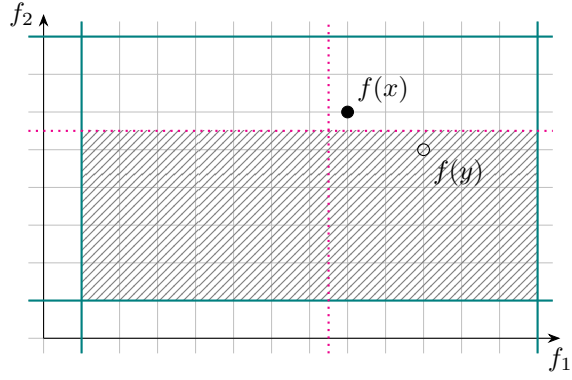

\begin{example} \label{ex:dpa}
    Consider a biobjective problem
    \begin{align*}
        \min_{x \in X \subseteq\{0,1\}^3} \quad & \begin{pmatrix}
            1721191 x_1 + 2417815 x_2 + 7529420 x_3 \\ 
            4855265 x_1 + 3892271 x_2 + 5357625 x_3
        \end{pmatrix} 
    \end{align*}
    The following two solutions and corresponding points are generated by DPA
    \begin{align*}
        &x = (0, 1, 0)  & f(x) = (2417815, 3892271)\\
&y = (-6.70085 \cdot 10^{-7}, 1, 4.18803 \cdot 10^{-7}) & f(y) = (2417817,3892270).
    \end{align*}
    The solution vectors are essentially the same, while the corresponding points differ and are not dominated by each other.
    Figure~\ref{fig:dpa} depicts the problematic behavior. When $x$ is generated, $f(x)$ is used to define boxes in the objective space that could potentially contain further non-dominated points. Since all points are assumed to be integer, the bounds of the boxes are derived from $f(x)$ but are decreased by $0.5$ to cut off $f(x)$ from the new boxes. However, the underlying solver now finds a new non-dominated point $f(y)$ in the hatched box, which is not a feasible point, since the entries of the solution vector $y$ are only considered binary up to the solver's accuracy.
    
\end{example}

The described problem can not be helped by rounding the newly found solution vector to a binary vector during the DPA, since this would also lead to re-generating the previously found non-dominated point.

\subsubsection{Instances and method}

We consider $15$ multi-objective instances\footnote{The instances and the implementation of the cutting plane algorithm can be found under: \url{https://git.zib.de/sriedmueller/objective-space-contraction}} for the Assignment-Problem (AP) and the Knapsack-Problem (KP) with $3$ to $5$ objective functions, see Table~\ref{tab:kpap}. 
All considered instances include contractable objectives and exhibit instability when DPA is applied. 
By instability, we refer to the unnecessary repeated generation of points that have already been found: in these cases, the number of duplicate points produced by DPA exceeds the number of actual unique solutions by at least a factor of 10.

\begin{table}[ht]
    \centering
    \begin{tabular}{lrrr}
     \hline
     name    &  p & \# variables & \# constraints\\
     \hline
    \texttt{AP\_p3\_n16\_v4} & 3 & 16 & 8 \\
    \texttt{AP\_p4\_n16\_v5} & 4 & 16 & 8 \\
    \texttt{AP\_p4\_n16\_v12} & 4 & 16 & 8 \\
    \texttt{AP\_p5\_n16\_v6} & 5 & 16 & 8 \\
    \texttt{AP\_p5\_n16\_v9} & 5 & 16 & 8 \\
     \hline
    \texttt{KP\_p4\_n10\_v0} & 4 & 10 & 1 \\
    \texttt{KP\_p4\_n10\_v9} & 4 & 10 & 1 \\
    \texttt{KP\_p4\_n15\_v3} & 4 & 15 & 1 \\
    \texttt{KP\_p4\_n15\_v7} & 4 & 15 & 1 \\
    \texttt{KP\_p5\_n10\_v0} & 5 & 10 & 1 \\
    \texttt{KP\_p5\_n10\_v1} & 5 & 10 & 1 \\
    \texttt{KP\_p5\_n15\_v0} & 5 & 15 & 1 \\
    \texttt{KP\_p5\_n15\_v3} & 5 & 15 & 1 \\
    \texttt{KP\_p5\_n30\_v1} & 5 & 30 & 1 \\
    \texttt{KP\_p5\_n30\_v2} & 5 & 30 & 1 \\
     \hline
    \end{tabular}
    \caption{Overview of the Instances for the Assignment and Knapsack Problem.}
    \label{tab:kpap}
\end{table}

As a base case, the presented instances in their original form are first solved by enumerating all possible solutions. While such an enumeration is only practicable for rather small instances, it provides a ground truth on the exact number of non-dominated points. The instances in their original form are then solved using an openly available DPA implementation.
We further solve the instances with DPA after preprocessing the objectives in the following ways:
\begin{itemize}
    \item The objectives are preprocessed by the scale-and-round heuristic from Section~\ref{sec:scaleandround}. We scale the objective coefficients by $10^k$ for $k=2,3,4$ and round up to avoid vanishing coefficients.
    \item The objectives are contracted as described in Section~\ref{sec:contraction}. For the runtime, we measure preprocessing time and solution time independently.
\end{itemize}
We compare the three versions with respect to runtime (preprocessing + solution time), the number of found solutions, and the quality of the found solutions.
A 3h time limit is set.

\subsubsection{Results}

We first evaluate the results regarding solution quantity and quality, and then compare the running and preprocessing times.

Table~\ref{tab:num_solutions} shows the number of found non-dominated points and how many of their corresponding solution vectors are actually distinct when rounded to exact binary vectors (recall Example~\ref{ex:dpa}). While the number of found points for DPA applied to the original problem is for the considered instances very large, the number of distinct solutions is comparably small (a difference in size of up to $7$ orders of magnitude). DPA is also not able to find all non-dominated points, since all instances but one (\texttt{KP\_p4\_n10\_v9}) do not terminate within the time limit.
In all instances, preprocessing by contraction or scaling before applying DPA terminates within the given time limit and finds more distinct solution vectors and corresponding non-dominated points than no preprocessing. For all instances, preprocessing by contraction leads to the largest number of distinct solutions. In 13 cases, the contraction allows DPA to find the exact number of non-dominated points; in 7 cases, the redundancy in the generated solutions is eliminated completely. As expected, an aggressive scaling leads to a decreased inaccuracy. In six cases, appropriate scaling led to the exact determination of the Pareto front. However, it is not clear how to choose the appropriate scaling factor in advance. 
As described in Example~\ref{ex:scaling_rounding}, 
Instance \texttt{AP\_p5\_n16\_v9} also shows that the scale-and-round heuristic can even introduce wrong non-dominated points.
It is worth noting that the contraction does not guarantee complete removal of instability, as observed in the two largest Knapsack instances. The success of the instability reduction proves to be more efficient, if the contraction factor is larger, since the two largest instances that still show some instability are also the instances with the smallest average contraction factor; compare Table~\ref{tab:gamma}. However, the instability is strongly reduced.
\begin{table}[ht]
\centering
\begin{tabular}{l c c c c c c}
\hline
name & $\gamma_1$ & $\gamma_2$ & $\gamma_3$ & $\gamma_4$ & $\gamma_5$ & avg \\
\hline
\texttt{AP\_p3\_n16\_v4}  & 0.991 & 0.991 & 0.990 & --    & --    & 0.990 \\
\texttt{AP\_p4\_n16\_v12} & 0.662 & 0.991 & 0.806 & 0.806 & --    & 0.816 \\
\texttt{AP\_p4\_n16\_v5}  & 0.991 & 0.965 & 0.965 & 0.965 & --    & 0.972 \\
\texttt{AP\_p5\_n16\_v6}  & 0.990 & 0.957 & 0.950 & 0.957 & 0.995 & 0.970 \\
\texttt{AP\_p5\_n16\_v9}  & 0.965 & 0.957 & 0.995 & 0.990 & 0.990 & 0.979 \\
\texttt{KP\_p4\_n10\_v0}  & 0.999 & 0.972 & 0.972 & 0.999 & --    & 0.986 \\
\texttt{KP\_p4\_n10\_v9}  & 0.999 & 0.998 & 0.999 & 0.999 & --    & 0.999 \\
\texttt{KP\_p4\_n15\_v3}  & 0.976 & 0.976 & 0.993 & 0.000 & --    & 0.736 \\
\texttt{KP\_p4\_n15\_v7}  & 0.993 & 0.978 & 0.978 & 0.993 & --    & 0.986 \\
\texttt{KP\_p5\_n10\_v0}  & 0.972 & 0.998 & 0.999 & 0.972 & 0.998 & 0.988 \\
\texttt{KP\_p5\_n10\_v1}  & 0.844 & 0.999 & 0.999 & 0.976 & 0.976 & 0.959 \\
\texttt{KP\_p5\_n15\_v0}  & 0.993 & 0.514 & 0.962 & 0.896 & 0.896 & 0.852 \\
\texttt{KP\_p5\_n15\_v3}  & 0.971 & 0.995 & 0.990 & 0.990 & 0.995 & 0.988 \\
\texttt{KP\_p5\_n30\_v1}  & 0.506 & 0.506 & 0.506 & 0.506 & 0.506 & 0.506 \\
\texttt{KP\_p5\_n30\_v2}  & 0.506 & 0.506 & 0.506 & 0.506 & 0.506 & 0.506 \\
\hline
\end{tabular}
\caption{Contraction factor $\gamma_i(f_i, g_i)$ and average per instance.}
\label{tab:gamma}
\end{table}

Table~\ref{tab:runtime} depicts the corresponding runtimes. The runtime for preprocessing by contraction is split into the actual solution time and the preprocessing time.  
All preprocessing techniques significantly reduced runtime, as the DPA applied to the original problem exceeded the time limit in all but one case. After preprocessing, the solution time after contraction exceeded the scaled variants in 2 instances; in the other instances, it performed competitively.
The preprocessing time for the contracted version is notably large in 4 instances. We observed that those are exactly the instances with the smallest average contraction factor.
As expected, an aggressive scaling leads to shorter run times.

\begin{table}[ht]
    \centering
\begin{tabular}{lrr|rr|rr|rr|rr|r}
\hline
 & \multicolumn{2}{r}{original} & \multicolumn{2}{r}{contraction} & \multicolumn{2}{r}{scale 2} & \multicolumn{2}{r}{scale 3} & \multicolumn{2}{r}{scale 4} & enum \\
instance & all & dist & all & dist & all & dist & all & dist & all & dist  \\
\hline
\texttt{AP\_p3\_n16\_v4} & 12291298 & 2 &
\textbf{7} & \textbf{7} &
\textbf{7} & \textbf{7} &
\textbf{7} & \textbf{7} &
3 & 3 & 7 \\

\texttt{AP\_p4\_n16\_v12} & 37536456 & \textbf{11} &
13 & \textbf{11} &
10 & 10 & 9 & 9 & 8 & 8 & 11 \\

\texttt{AP\_p4\_n16\_v5} & 10415768 & 1 &
18 & \textbf{16} &
\textbf{16} & \textbf{16} &
\textbf{16} & \textbf{16} &
\textbf{16} & \textbf{16} & 16 \\

\texttt{AP\_p5\_n16\_v6} & 12015465 & 1 &
19 & \textbf{14} &
\textbf{14} & \textbf{14} &
\textbf{14} & \textbf{14} &
\textbf{14} & \textbf{14} & 14 \\

\texttt{AP\_p5\_n16\_v9} & 9643184 & 10 &
\textbf{20} & \textbf{20} &
\textbf{20} & \textbf{20} &
18 & 18 &
12 & 12 & 20 \\

\texttt{KP\_p4\_n10\_v0} & 4251844 & 3 &
\textbf{49} & \textbf{49} &
50 & 50 & 42 & 42 & 18 & 18 & 49 \\

\texttt{KP\_p4\_n10\_v9} & 1100 & \textbf{23} &
\textbf{23} & \textbf{23} &
\textbf{23} & \textbf{23} &
\textbf{23} & \textbf{23} &
\textbf{23} & \textbf{23} & 23 \\

\texttt{KP\_p4\_n15\_v3} & 15596776 & 5 &
170 & \textbf{169} &
142 & 142 & 145 & 145 & 137 & 137 & 169 \\

\texttt{KP\_p4\_n15\_v7} & 10928060 & 10 &
130 & \textbf{111} &
96 & 96 & 93 & 93 & 90 & 90 & 111 \\

\texttt{KP\_p5\_n10\_v0} & 12068799 & 9 &
\textbf{51} & \textbf{51} &
49 & 48 & 43 & 43 & 30 & 30 & 51 \\

\texttt{KP\_p5\_n10\_v1} & 16870007 & 8 &
\textbf{29} & \textbf{29} &
31 & \textbf{29} &
23 & 23 & 8 & 8 & 29 \\

\texttt{KP\_p5\_n15\_v0} & 12748705 & 11 &
550 & \textbf{421} &
413 & 413 & 406 & 406 & 182 & 182 & 421 \\

\texttt{KP\_p5\_n15\_v3} & 16319290 & 5 &
\textbf{183} & \textbf{183} &
171 & 171 & 145 & 145 & 113 & 113 & 183 \\

\texttt{KP\_p5\_n30\_v1} & 9896195 & 11 &
5982 & \textbf{4603} &
3373 & 3373 &
1075 & 1075 &
58 & 58 & 4638 \\

\texttt{KP\_p5\_n30\_v2} & 16544856 & 2 &
9441 & \textbf{7320} &
6105 & 6105 &
2907 & 2907 &
50 & 50 & 7326 \\
\hline
\end{tabular}

    \caption{Depicted are the number of found points per method (all) and how many of the corresponding solutions are actually distinct (dist), when the solution vectors have been rounded to binary values. Emphasized are those results that are closest to the ground truth based on enumeration (enum).}
    \label{tab:num_solutions}
\end{table}

\begin{table}[ht]
\centering
\begin{tabular}{lr|rr|rrr}
\hline
\multicolumn{2}{c}{} & \multicolumn{2}{c}{contraction} & \multicolumn{3}{c}{scale-and-round}  \\
instance & original & solution & preprocessing & scale 2 & scale 3 & scale 4 \\
\hline
\texttt{AP\_p3\_n16\_v4} & limit & 0.009 & 0.338 & 0.009 & 0.008 & 0.007 \\
\texttt{AP\_p4\_n16\_v12} & limit & 0.028 & 1109.200 & 0.022 & 0.017 & 0.013 \\
\texttt{AP\_p4\_n16\_v5} & limit & 0.050 & 0.995 & 0.052 & 0.047 & 0.044 \\
\texttt{AP\_p5\_n16\_v6} & limit & 0.078 & 0.956 & 0.071 & 0.047 & 0.038 \\
\texttt{AP\_p5\_n16\_v9} & limit & 0.137 & 0.775 & 0.101 & 0.097 & 0.039 \\
\texttt{KP\_p4\_n10\_v0} & limit & 0.414 & 0.686 & 0.380 & 0.151 & 0.062 \\
\texttt{KP\_p4\_n10\_v9} & 2.766 & 0.302 & 1.720 & 0.284 & 0.292 & 0.251 \\
\texttt{KP\_p4\_n15\_v3} & limit& 2.014 & 600.459 & 1.711 & 1.576 & 1.144 \\
\texttt{KP\_p4\_n15\_v7} & limit & 1.255 & 0.566 & 0.937 & 0.713 & 0.684 \\
\texttt{KP\_p5\_n10\_v0} & limit & 0.733 & 0.643 & 0.635 & 0.335 & 0.160 \\
\texttt{KP\_p5\_n10\_v1} & limit & 0.289 & 0.947 & 0.236 & 0.137 & 0.032 \\
\texttt{KP\_p5\_n15\_v0} & limit & 29.480 & 0.749 & 21.781 & 17.307 & 5.164 \\
\texttt{KP\_p5\_n15\_v3} & limit & 6.002 & 1.262 & 3.690 & 2.619 & 1.525 \\
\texttt{KP\_p5\_n30\_v1} & limit & 493.340 & 825.472 & 204.745 & 53.562 & 0.969 \\
\texttt{KP\_p5\_n30\_v2} & limit & 1020.740 & 817.153 & 577.013 & 186.239 & 1.284 \\
\hline
\end{tabular}
\caption{Solution times [s] for DPA on the original problem and DPA applied to a preprocessed version.}
\label{tab:runtime}
\end{table}

\section{Conclusion}

We introduced a new problem class of finding the smallest possible objective coefficients for a \ref{eq:moco} problem.
To solve the introduced problem, we presented an exact transformation approach that can be formulated as an integer linear program and solved using a cutting-plane algorithm. In mixed-integer programming, the contraction can still be applied to those objectives that involve only binary variables.

In a first computational study, the possibility of contracting an objective function has been assessed. For variants of randomly generated objective coefficients, the contractability of the objectives decreases, and the runtime of the contraction increases with the number of coefficients. The objectives were more likely to be contractable if the randomly generated coefficients showed a wide spread in order of magnitude.
Since applied settings often involve a large number of objective coefficients, the presented contraction method is likely applicable only in rare cases.

In a second study, we investigated unstable instances of multi-objective combinatorial problems. Heuristic scaling of the objective coefficients improved the stability of DPA immensely, but comes with an decrease in inaccuracy. The exact contraction reduced the instability effectively and helped DPA to find the largest number of distinct solutions.

In total, the presented contraction approach can be applied if the involved objective functions are contractable. In this case, the stability of the DPA algorithm can be increased significantly.
If the exactness of the solution set is not important, the heuristic scaling approach can be an efficient option.

The presented preprocessing technique should, of course, be considered only as one of several preprocessing or modeling decisions to increase efficiency. We generally recommend assessing whether the precision of the objective coefficients is truly necessary or can be relaxed. In particular, when dealing with a wide range of coefficient values, the relevance of enforcing integrality on the corresponding variables should be reconsidered.

For future research, it would be of interest to predict more accurately which properties an objective function must satisfy to be contractable. The cutting-plane algorithm could potentially be improved by deriving a tighter dual bound. In general, it remains an open question whether contracting the objective function in multi-objective settings provides benefits beyond stabilizing the DPA solution process.

\section{Acknowledgement}
\acknowledgments

\bibliographystyle{apalike}
\bibliography{references}

@article{NeusselStein2025,
  author  = {F. Neussel and O. Stein},
  title   = {On image space transformations in multiobjective optimization},
  journal = {Optimization Online},
  year    = {2025},
  note    = {Preprint},
  url     = {https://optimization-online.org/?p=30125}
}

@article{Daechert2024,
  title={A simple, efficient and versatile objective space algorithm for multiobjective integer programming},
  author={Dächert, Kerstin and Fleuren, Tino and Klamroth, Kathrin},
  journal={Mathematical Methods of Operations Research},
  volume={100},
  pages={351–384},
  year={2024},
  url={https://doi.org/10.1007/s00186-023-00841-0},
}

@article{Ozlen2014,
  title={Multi-Objective Integer Programming: An Improved Recursive Algorithm},
  author={Özlen, Melih and Burton, Benjamin A. and MacRae, Cameron A. G.},
  journal={Journal of Optimization Theory and Applications},
  volume={160},
  number={2},
  pages={470–482},
  year={2014},
  url={https://doi.org/10.1007/s10957-013-0364-y}
}

@Article{Kirlik2014,
  author={Kirlik, Gokhan and Sayın, Serpil},
  title={{A new algorithm for generating all nondominated solutions of multiobjective discrete optimization problems}},
  journal={European Journal of Operational Research},
  year={2014},
  volume={232},
  number={3},
  pages={479-488},
  month={},
  doi={10.1016/j.ejor.2013.08.00},
  url={https://ideas.repec.org/a/eee/ejores/v232y2014i3p479-488.html}
}

@Article{Boland2017,
  author={Boland, Natashia and Charkhgard, Hadi and Savelsbergh, Martin},
  title={{A new method for optimizing a linear function over the efficient set of a multiobjective integer program}},
  journal={European Journal of Operational Research},
  year={2017},
  volume={260},
  number={3},
  pages={904-919},
  month={},
  doi={10.1016/j.ejor.2016.02.03},
  url={https://ideas.repec.org/a/eee/ejores/v260y2017i3p904-919.html}
}

@misc{Bauss2023,
      title={Adaptive Improvements of Multi-Objective Branch and Bound}, 
      author={Julius Bauß and Sophie N. Parragh and Michael Stiglmayr},
      year={2023},
      eprint={2312.12192},
      archivePrefix={arXiv},
      primaryClass={math.OC},
      url={https://arxiv.org/abs/2312.12192}, 
}

@book{Ehrgott2005,
    author = {Matthias Ehrgott},
    title = {Multicriteria Optimization},
    publisher = {Springer Berlin, Heidelberg},
    year = {2005},
    DOI = {https://doi.org/10.1007/3-540-27659-9}
}

@article{Figueira2016,
author = {Figueira, José and Fonseca, Carlos and Halffmann, Pascal and Klamroth, Kathrin and Paquete, Luis and Ruzika, Stefan and Schulze, Britta and Stiglmayr, Michael and Willems, David},
year = {2016},
month = {01},
pages = {},
title = {Easy to say they are Hard, but Hard to see they are Easy- Towards a Categorization of Tractable Multiobjective Combinatorial Optimization Problems},
volume = {24},
journal = {Journal of Multi-Criteria Decision Analysis},
doi = {10.1002/mcda.1574}
}

@article{Allmendinger2022,
title = {What if we increase the number of objectives? Theoretical and empirical implications for many-objective combinatorial optimization},
journal = {Computers \& Operations Research},
volume = {145},
pages = {105857},
year = {2022},
issn = {0305-0548},
doi = {https://doi.org/10.1016/j.cor.2022.105857},
url = {https://www.sciencedirect.com/science/article/pii/S0305054822001319},
author = {Richard Allmendinger and Andrzej Jaszkiewicz and Arnaud Liefooghe and Christiane Tammer}
}

@manual{ibm_cplex,
  title        = {IBM ILOG CPLEX Optimization Studio},
  author       = {{IBM}},
  year         = {2019},
  note         = {Version 12.10.0},
  url          = {https://www.ibm.com/products/ilog-cplex-optimization-studio}
}

@article{klamroth2023,
title = {Ordinal optimization through multi-objective reformulation},
journal = {European Journal of Operational Research},
volume = {311},
number = {2},
pages = {427-443},
year = {2023},
issn = {0377-2217},
doi = {https://doi.org/10.1016/j.ejor.2023.04.042},
url = {https://www.sciencedirect.com/science/article/pii/S0377221723003399},
author = {Kathrin Klamroth and Michael Stiglmayr and Julia Sudhoff}
}

@misc{loehne2025,
      title={Low-Rank Multi-Objective Linear Programming}, 
      author={Andreas Löhne and Pascal Zillmann},
      year={2025},
      eprint={2508.20880},
      archivePrefix={arXiv},
      primaryClass={math.OC},
      url={https://arxiv.org/abs/2508.20880}, 
}

@article{Boekler2017,
author = {Bökler, Fritz and Ehrgott, Matthias and Morris, Christopher and Mutzel, Petra},
title = {Output-sensitive complexity of multiobjective combinatorial optimization},
journal = {Journal of Multi-Criteria Decision Analysis},
volume = {24},
number = {1-2},
pages = {25-36},
doi = {https://doi.org/10.1002/mcda.1603},
url = {https://onlinelibrary.wiley.com/doi/abs/10.1002/mcda.1603},
eprint = {https://onlinelibrary.wiley.com/doi/pdf/10.1002/mcda.1603},
year = {2017}
}

@mastersthesis{Kof2024,
    author = {Gökhan Kof},
    title = {Eliminating Objective Functions and Warm Starting ALgorithms Using Projections in Multi-objective Optimization},
    school = {Ko\c{c} University},
    year = {2024}
}

@article{brockhoff2009,
  author       = {Brockhoff, D. and Zitzler, E.},
  title        = {Objective reduction in evolutionary multiobjective optimization: theory and applications},
  journal      = {Evolutionary Computation},
  year         = {2009},
  volume       = {17},
  number       = {2},
  pages        = {135--166},
  doi          = {10.1162/evco.2009.17.2.135},
}

@article{VAZQUEZ2018382,
title = {MILP method for objective reduction in multi-objective optimization},
journal = {Computers \& Chemical Engineering},
volume = {108},
pages = {382-394},
year = {2018},
issn = {0098-1354},
doi = {https://doi.org/10.1016/j.compchemeng.2017.10.021},
url = {https://www.sciencedirect.com/science/article/pii/S0098135417303769},
author = {Daniel Vázquez and María J. Fernández-Torres and Rubén Ruiz-Femenia and Laureano Jiménez and José A. Caballero}
}

@INPROCEEDINGS{Cheung2014,
  author={Cheung, Yiu-ming and Gu, Fangqing},
  booktitle={2014 IEEE Congress on Evolutionary Computation (CEC)}, 
  title={Online objective reduction for many-objective optimization problems}, 
  year={2014},
  volume={},
  number={},
  pages={1165-1171},
  doi={10.1109/CEC.2014.6900548}}

@article{Li1996,
  author       = {Li, D.},
  title        = {Convexification of a noninferior frontier},
  journal      = {Journal of Optimization Theory and Applications},
  year         = {1996},
  volume       = {88},
  number       = {1},
  pages        = {177--196},
  doi          = {10.1007/BF02192028},
  url          = {https://doi.org/10.1007/BF02192028},
  issn         = {1573-2878}
}

@article{Li1998,
  author       = {Li, D. and Biswal, M. P.},
  title        = {Exponential Transformation in Convexifying a Noninferior Frontier and Exponential Generating Method},
  journal      = {Journal of Optimization Theory and Applications},
  year         = {1998},
  volume       = {99},
  number       = {1},
  pages        = {183--199},
  doi          = {10.1023/A:1021708412776},
  url          = {https://doi.org/10.1023/A:1021708412776},
  issn         = {1573-2878}
}

@article{Zarepisheh2017,
  author       = {Zarepisheh, Masoud and Pardalos, Panos M.},
  title        = {An equivalent transformation of multi-objective optimization problems},
  journal      = {Annals of Operations Research},
  year         = {2017},
  volume       = {249},
  number       = {1},
  pages        = {5--15},
  doi          = {10.1007/s10479-014-1782-4},
  url          = {https://doi.org/10.1007/s10479-014-1782-4},
  issn         = {1572-9338}
}

@article{Romeijn2004,
  author       = {Romeijn, H. E. and Dempsey, J. F. and Li, J. G.},
  title        = {A unifying framework for multi-criteria fluence map optimization models},
  journal      = {Physics in Medicine and Biology},
  year         = {2004},
  month        = {May 21},
  volume       = {49},
  number       = {10},
  pages        = {1991--2013},
  doi          = {10.1088/0031-9155/49/10/011},
  pmid         = {15214537}
}

@article{Hirschberger2005,
author = {Hirschberger, M.},
year = {2005},
month = {06},
pages = {283-304},
title = {Connectedness of efficient points in convex and convex transformable vector optimization},
volume = {54},
journal = {Optimization},
doi = {10.1080/02331930500096270}
}

@article{Marler2005,
author = {Marler, R. and Arora, Jasbir},
year = {2005},
month = {09},
pages = {551-570},
title = {Function-transformation methods for multi-objective Optimization},
volume = {37},
journal = {Engineering Optimization},
doi = {10.1080/03052150500114289}
}

@article{BOLAND2019858,
title = {Preprocessing and cut generation techniques for multi-objective binary programming},
journal = {European Journal of Operational Research},
volume = {274},
number = {3},
pages = {858-875},
year = {2019},
issn = {0377-2217},
doi = {https://doi.org/10.1016/j.ejor.2018.10.034},
url = {https://www.sciencedirect.com/science/article/pii/S0377221718308877},
author = {Natashia Boland and Hadi Charkhgard and Martin Savelsbergh}
}

@InProceedings{Jabs2023,
  author =	{Jabs, Christoph and Berg, Jeremias and Ihalainen, Hannes and J\"{a}rvisalo, Matti},
  title =	{{Preprocessing in SAT-Based Multi-Objective Combinatorial Optimization}},
  booktitle =	{29th International Conference on Principles and Practice of Constraint Programming (CP 2023)},
  pages =	{18:1--18:20},
  series =	{Leibniz International Proceedings in Informatics (LIPIcs)},
  ISBN =	{978-3-95977-300-3},
  ISSN =	{1868-8969},
  year =	{2023},
  volume =	{280},
  editor =	{Yap, Roland H. C.},
  publisher =	{Schloss Dagstuhl -- Leibniz-Zentrum f{\"u}r Informatik},
  address =	{Dagstuhl, Germany},
  URL =		{https://drops.dagstuhl.de/entities/document/10.4230/LIPIcs.CP.2023.18},
  URN =		{urn:nbn:de:0030-drops-190553},
  doi =		{10.4230/LIPIcs.CP.2023.18}
}

@article{Savelsbergh1994,
  author    = {Savelsbergh, M. W. P.},
  title     = {Preprocessing and Probing Techniques for Mixed Integer Programming Problems},
  journal   = {ORSA Journal on Computing},
  year      = {1994},
  volume    = {6},
  number    = {4},
  pages     = {445--454},
  doi       = {10.1287/ijoc.6.4.445}
}

@incollection{FUGENSCHUH200569,
title = {Computational Integer Programming and Cutting Planes},
editor = {K. Aardal and G.L. Nemhauser and R. Weismantel},
series = {Handbooks in Operations Research and Management Science},
publisher = {Elsevier},
volume = {12},
pages = {69-121},
year = {2005},
booktitle = {Discrete Optimization},
issn = {0927-0507},
doi = {https://doi.org/10.1016/S0927-0507(05)12002-7},
url = {https://www.sciencedirect.com/science/article/pii/S0927050705120027},
author = {Armin Fügenschuh and Alexander Martin}
}

@book{chen2011applied,
  title={Applied integer programming: modeling and solution},
  author={Chen, Der-San and Batson, Robert G and Dang, Yu},
  year={2011},
  publisher={John Wiley \& Sons}
}

@article{GAL1977176,
title = {Redundant objective functions in linear vector maximum problems and their determination},
journal = {European Journal of Operational Research},
volume = {1},
number = {3},
pages = {176-184},
year = {1977},
issn = {0377-2217},
doi = {https://doi.org/10.1016/0377-2217(77)90025-X},
url = {https://www.sciencedirect.com/science/article/pii/037722177790025X},
author = {Tomas Gal and Heiner Leberling}
}

@article{LINDROTH20101519,
title = {Approximating the Pareto optimal set using a reduced set of objective functions},
journal = {European Journal of Operational Research},
volume = {207},
number = {3},
pages = {1519-1534},
year = {2010},
issn = {0377-2217},
doi = {https://doi.org/10.1016/j.ejor.2010.07.004},
url = {https://www.sciencedirect.com/science/article/pii/S0377221710004868},
author = {Peter Lindroth and Michael Patriksson and Ann-Brith Strömberg}
}

@article{Malinowska2008,
  author    = {Malinowska, A. B. and Torres, D. F. M.},
  title     = {Computational Approach to Essential and Nonessential Objective Functions in Linear Multicriteria Optimization},
  journal   = {Journal of Optimization Theory and Applications},
  year      = {2008},
  volume    = {139},
  number    = {3},
  pages     = {577--590},
  doi       = {10.1007/s10957-008-9397-z}
}

@article{Thoai2012,
  author    = {Thoai, Nguyen Van},
  title     = {Criteria and Dimension Reduction of Linear Multiple Criteria Optimization Problems},
  journal   = {Journal of Global Optimization},
  year      = {2012},
  volume    = {52},
  number    = {3},
  pages     = {499--508},
  doi       = {10.1007/s10898-011-9764-4}
}

@article{Halffmann2022,
author = {Halffmann, Pascal and Schäfer, Luca E. and Dächert, Kerstin and Klamroth, Kathrin and Ruzika, Stefan},
title = {Exact algorithms for multiobjective linear optimization problems with integer variables: A state of the art survey},
journal = {Journal of Multi-Criteria Decision Analysis},
volume = {29},
number = {5-6},
pages = {341-363},
doi = {https://doi.org/10.1002/mcda.1780},
url = {https://onlinelibrary.wiley.com/doi/abs/10.1002/mcda.1780},
eprint = {https://onlinelibrary.wiley.com/doi/pdf/10.1002/mcda.1780},
year = {2022}
}

@article{Antunes2024,
  author    = {Antunes, C. H. and Fonseca, C. M. and Paquete, L. and others},
  title     = {Special Issue on Exact and Approximation Methods for Mixed-Integer Multi-Objective Optimization},
  journal   = {Mathematical Methods of Operations Research},
  year      = {2024},
  volume    = {100},
  number    = {1},
  pages     = {1--4},
  doi       = {10.1007/s00186-024-00874-z}
}

@INPROCEEDINGS{Riquelme2015,
  author={Riquelme, Nery and Von Lücken, Christian and Baran, Benjamin},
  booktitle={2015 Latin American Computing Conference (CLEI)}, 
  title={Performance metrics in multi-objective optimization}, 
  year={2015},
  volume={},
  number={},
  pages={1-11},
  doi={10.1109/CLEI.2015.7360024}}

@article{Helfrich2024,
  author    = {Stefan Helfrich and Anne Herzel and Sebastian Ruzika and Clemens Thielen},
  title     = {Using scalarizations for the approximation of multiobjective optimization problems: towards a general theory},
  journal   = {Mathematical Methods of Operations Research},
  year      = {2024},
  volume    = {100},
  pages     = {27--63},
  doi       = {10.1007/s00186-023-00823-2},
  url       = {https://doi.org/10.1007/s00186-023-00823-2}
}

@book{Higham2002,
author = {Higham, Nicholas J.},
title = {Accuracy and Stability of Numerical Algorithms},
publisher = {Society for Industrial and Applied Mathematics},
year = {2002},
doi = {10.1137/1.9780898718027},
address = {},
edition   = {Second},
URL = {https://epubs.siam.org/doi/abs/10.1137/1.9780898718027},
eprint = {https://epubs.siam.org/doi/pdf/10.1137/1.9780898718027}
}

@book{Widrow_Kollár_2008, place={Cambridge}, title={Quantization Noise: Roundoff Error in Digital Computation, Signal Processing, Control, and Communications}, publisher={Cambridge University Press}, author={Widrow, Bernard and Kollár, István}, year={2008}}

@techreport{TenfeldePodehl2003,
  author       = {Tenfelde-Podehl, D.},
  title        = {A Recursive Algorithm for Multiobjective Combinatorial Optimization Problems with {Q} Criteria},
  institution  = {Institut für Mathematik, Technische Universität Graz},
  year         = {2003},
  type         = {Technical Report}
}

@article{OZLEN200925,
title = {Multi-objective integer programming: A general approach for generating all non-dominated solutions},
journal = {European Journal of Operational Research},
volume = {199},
number = {1},
pages = {25-35},
year = {2009},
issn = {0377-2217},
doi = {https://doi.org/10.1016/j.ejor.2008.10.023},
url = {https://www.sciencedirect.com/science/article/pii/S0377221708009624},
author = {Melih Özlen and Meral Azizoğlu}
}

@article{LAUMANNS2006932,
title = {An efficient, adaptive parameter variation scheme for metaheuristics based on the epsilon-constraint method},
journal = {European Journal of Operational Research},
volume = {169},
number = {3},
pages = {932-942},
year = {2006},
issn = {0377-2217},
doi = {https://doi.org/10.1016/j.ejor.2004.08.029},
url = {https://www.sciencedirect.com/science/article/pii/S0377221704005715},
author = {Marco Laumanns and Lothar Thiele and Eckart Zitzler}
}

@article{Ehrgott2006,
  author    = {Ehrgott, M.},
  title     = {A Discussion of Scalarization Techniques for Multiple Objective Integer Programming},
  journal   = {Annals of Operations Research},
  volume    = {147},
  pages     = {343--360},
  year      = {2006},
  doi       = {10.1007/s10479-006-0074-z},
  url       = {https://doi.org/10.1007/s10479-006-0074-z}
}

\end{document}